\documentclass{amsart}
\usepackage{subfigure, xypic}

\usepackage{amssymb, 
     amsmath, 
     amsthm,
     graphicx,
     amscd,
     multind,
     eufrak
     }

\theoremstyle{plain}
\newtheorem{thm}{Theorem}[section]
\newtheorem{lm}[thm]{Lemma}
\newtheorem{cor}[thm]{Corollary}
\newtheorem{prop}[thm]{Proposition}

\newtheorem{prb}[thm]{Problem}

\theoremstyle{definition}
\newtheorem{de}[thm]{Definition}
\newtheorem{ex}[thm]{Example}

\newtheorem{re}[thm]{Remark}

\newcommand{\CC}{{\mathbb C}}
\newcommand{\RR}{{\mathbb R}}

\newcommand{\ZZ}{{\mathbb Z}}
\newcommand{\NN}{{\mathbb N}}

\newcommand{\PP}{{\mathbb P}}

\newcommand{\Rb}{\overline{\RR}}

\newcommand{\im}{\operatorname{im}}

	{\begin{figure} \begin{center}}%
	{\end{center} \end{figure}}

\newcommand{\la}{\langle}
\newcommand{\ra}{\rangle}

\newcommand{\Wedge}{\bigwedge\nolimits}
\newcommand{\wt}{\mathrm{wt}}
\newcommand{\ini}{\mathrm{in}}
\newcommand{\lx}{{\bar{x}}}
\newcommand{\ly}{{\bar{y}}}

\newcommand{\liea}[1]{\mathfrak{#1}}
\newcommand{\lieg}[1]{\mathrm{#1}}

\newcommand{\Aff}{\operatorname{Aff}\nolimits}
\newcommand{\LP}{\text{\sc LinearPartition}}
\newcommand{\AP}{\text{\sc AffinePartition}}
\newcommand{\VP}{\text{\sc VoronoiPartition}}
\newcommand{\Vor}{\operatorname{Vor}\nolimits}

\begin{document}
\title{A tropical approach to secant dimensions}
\author{Jan Draisma}
\address{
Department of Mathematics and Computer Science\\
Technische Universiteit Eindhoven\\
P.O. Box 513, 5600 MB Eindhoven\\
The Netherlands}
\email{j.draisma@tue.nl}

\maketitle

\section*{Abstract}
Tropical geometry yields good lower bounds, in terms of certain
combinatorial-polyhedral optimisation problems, on the dimensions
of secant varieties. In particular, it gives an attractive pictorial
proof of the theorem of Hirschowitz that all Veronese embeddings of
the projective plane except for the quadratic one and the quartic
one are non-defective; this proof might be generalisable to cover
all Veronese embeddings, whose secant dimensions are known from the
ground-breaking but difficult work of Alexander and Hirschowitz. Also,
the non-defectiveness of certain Segre embeddings is proved, which
cannot be proved with the rook covering argument already known in the
literature. Short self-contained introductions to secant varieties and
the required tropical geometry are included.

\section*{Acknowledgments}
I thank Tony Geramita, Rick Miranda, Bernd Sturmfels, and Seth Sullivant
very much for their responses to an earlier version of this paper. In
particular, they pointed out several very nice papers on secant varieties
of which I was unaware. Also, their responses raise the question
of how the tropical approach relates to the degeneration approach
of Sturmfels-Sullivant and Dumitrescu-Miranda. This question and the
combinatorial-polyhedral problems introduced below will certainly lead
to further research. Finally, I thank Hannah Markwig for pointing out
a gap in the proof of Theorem \ref{thm:Characterisations} in an earlier
version of this paper, and Marius van der Put who pointed out a theorem
about affinoid algebras which can be used to fill that gap.

\section{Introduction} \label{sec:Introduction}
Secant varieties are rather classical objects of study in algebraic
geometry: given a closed subvariety $X$ of some projective space $\PP^m$,
and given a natural number $k$, one tries to describe the union of all
subspaces of $\PP^m$ that are spanned by $k$ points on $X$. We call the
Zariski closure of this union the $k$-th {\em secant variety} of $X$,
and denote it by $kX$. To avoid confusion: some authors call this the
$(k-1)$st secant variety. So in this paper $2X$ is the variety of secant
lines, traditionally called {\em the} secant variety of $X$. We will
refer to all $kX$ as (higher) secant varieties, and to their dimensions as
(higher) secant dimensions. The standard reference for secant varieties
is \cite{Zak93}.

Already the most basic of all questions about the secant varieties of $X$
poses unexpected challenges, namely: what are their dimensions? This
question is of particular interest when $X$ is a {\em minimal orbit}
in a representation space of a reductive group. These minimal orbits
comprise Segre embeddings of products of projective spaces, Pl\"ucker
embeddings of Grassmannians, and Veronese embeddings of projective spaces;
see Section \ref{sec:MinOrbit}. Among these instances, only the secant
dimensions of the Veronese embeddings are completely known, from the
ground-breaking work of Alexander and Hirschowitz \cite{Alexander88,
Alexander92, Alexander95, Hirschowitz85}. Secant dimensions of
Segre powers of the projective line are almost entirely known
\cite{Catalisano05b}.

This paper introduces a new approach to secant dimensions,
based on tropical geometry. Tropical geometry is well known
as a tool for transforming algebraic-geometric questions
into polyhedral-combinatorial ones. Recommended references are
\cite{Bieri84,Einsiedler04,Mikhalkin04,Pachter05,RichterGebert05,
Speyer04} and the references therein---however, all background in tropical
geometry needed here is reviewed in Section \ref{sec:Tropical}.

In Sections \ref{sec:Join} and \ref{sec:Secant} I present the tropical
lower bounds on secant dimensions in terms of certain polyhedral
optimisation problems. After a review of the necessary tropical
geometry in Section \ref{sec:Tropical} we prove the lower bounds in
Section \ref{sec:Proofs}.  In Section \ref{sec:MinOrbit} I recall the
notion of minimal orbits, and give two lower bounds on their secant
dimensions. One of them is well known in special cases; for instance,
it uses rook coverings for Segre varieties, and a variation of these for
Grassmannians \cite{Catalisano2002,Catalisano05,Ehrenborg00,Sturmfels05}.
The other seems to be good for Segre products of Veronese embeddings.

Then in Section \ref{sec:Applications} we apply the tropical lower bounds
to Segre varieties, Veronese embeddings, and Grassmannians, and show that
they are better than the bounds from Section \ref{sec:MinOrbit}.  As an
example, we reprove the theorem that all but two Veronese embeddings
of the projective plane are non-defective; this was proved earlier
by Hirschowitz \cite{Hirschowitz85} using his ``Horace method'' and
by Miranda and Dumitrescu using degenerations (private communication).
Also, I give a nice proof that the $6$-fold Segre power of the projective
line is non-defective; this cannot be proved using rook coverings alone,
and is the first case not covered by \cite{Catalisano05b}.  Finally, Seth
Sullivant and Bernd Sturmfels pointed out the paper \cite{Develin06}
to me, in which tropical secant varieties of ordinary linear spaces
are considered. The precise relation between Develin's tropical secant
varieties and the tropicalisation of secant varieties is still unclear
to me, though under reasonable conditions the former should be contained
in the latter.

In conclusion, the tropical approach is conceptually very simple, but
shows very promising results when tested on concrete examples. However,
it also raises many intriguing combinatorial-polyhedral optimisation
problems; I do not know of any efficient programs solving these.

\section{Joins, secant varieties, and first results}
\label{sec:Join}

Rather than projective varieties, we consider closed cones in affine
spaces. So let $K$ be an algebraically closed field of characteristic $0$,
let $V$ be a finite-dimensional vector space over $K$, and let $C,D$ be
closed cones: Zariski-closed subsets of $V$ that are closed under scalar
multiplication. Then we define the {\em join} of $C$ and $D$ as follows:
\[ C+D:=\overline{\{c+d \mid c \in C, d \in D\}}. \]
Note that in taking the closure we ignore the subtle question of which
elements of $C+D$ can actually be written as $c+d$ with $c \in C$
and $d \in D$; in this paper we are only interested in dimensions,
and hence there is no harm in taking the closure. 
There is an obvious upper bound on the dimension of $C+D$, namely
$\min\{\dim C + \dim D,\dim V\}$---indeed, the summation map $C \times
D \rightarrow C+D$ is dominant. We call this upper bound the {\em
expected dimension} of $C+D$. If $C+D$ has strictly lower dimension than
expected, then we call $C+D$ {\em defective}; otherwise, we call $C+D$
{\em non-defective}. The difference $\min\{\dim C+\dim D,\dim V\}-\dim(C+D)$
is called the {\em defect}.

Taking the join is an associative (and commutative) operation on
closed cones in $V$, so given $k$ closed cones $C_1,\ldots,C_k$, their join
$C_1+\ldots+C_k$ is well-defined. Again, we call this join defective
or non-defective according as its dimension is smaller than or equal to
$\min\{\dim V, \sum_i \dim C_i\}$. 

In particular, taking all $C_i$ equal to a single closed cone $C$ we
obtain $kC$, called the $k$th {\em secant variety} of $C$. The defect
of $kC$, also called the $k$-th secant defect of $C$, is defined in a
slightly different manner: it is the difference $\min \{\dim V,\dim C +
\dim (k-1)C\}-\dim kC$; hence if $k\dim C < \dim V$, then the difference
$k \dim C - \dim kC$ is the {\em sum} of all $l$-defects for $l \leq k$.
We call $kC$ defective if its defect is positive, and non-defective
otherwise. Finally, we call $C$ itself defective if and only if $kC$
is defective for some $k \geq 0$, and we call the numbers $\dim kC,\
k \in \NN$ the {\em secant dimensions} of $C$. The standard
reference for joins and secant varieties is \cite{Zak93}.

Typically, one considers a class of cones (e.g., the cones over
Grassmannians), one knows a short explicit list of defective secant
varieties of cones in this class, and wishes to prove that all other
secant varieties of cones in this class are non-defective. One then needs
{\em lower bounds} on secant dimensions that are in fact {\em equal}
to the expected dimensions---so that one can conclude equality.

Our approach towards such lower bounds focuses on the following, special
si\-tuation: suppose that $C_1,\ldots,C_k$ are closed cones in $V$,
and single out a basis $e_1,\ldots,e_n$ of $V$. The method depends on
this basis, but in our applications there will be natural bases to
work with.  Let $y_1,\ldots,y_n$ be the dual basis of $V^*$. Assume
for simplicity that none of the $C_i$ is contained in any coordinate
hyperplane $\{y_b=0\}$. Furthermore, suppose that for each $i$ we have
a finite-dimensional vector space $V_i$ over $K$, again with a fixed
basis $x=(x_1,\ldots,x_{m_i})$ of $V_i^*$, and a polynomial map $f_i:V_i
\rightarrow V$ that maps $V_i$ dominantly into $C_i$. In
particular, every $C_i$ is irreducible.

Write each $f_i$, relative to the bases of $V_i$ and $V$, as a list
$(f_{i,b})_{b=1}^n$ of polynomials $f_{i,b} \in K[x_1,\ldots,x_{m_i}]$;
the fact that we use the same letter $x$ to indicate coordinates on the
distinct $V_i$ will not lead to any confusion. For every $i=1,\ldots,k$
and $b=1,\ldots,n$, let $l_{i,b}$ be the piecewise linear function
$\RR^{m_i} \rightarrow \RR$ defined as follows: Write
\[ f_{i,b}=\sum_\alpha c_\alpha x^\alpha, \]
where the sum runs over all multi-indices $\alpha \in \NN^{m_i}$ for
which $c_\alpha$ is non-zero; we view these $\alpha$ as row vectors,
i.e., we will tacitly regard $\NN^{m_i}$ as a subset of $(\RR^{m_i})^*$.
Note that this is not an empty sum by the assumption that $C_i$ not lie
in $\{y_b=0\}$. Then $l_{i,b}$ is defined by
\[ l_{i,b}(v):=\min_\alpha \la v,\alpha \ra, \quad v \in \RR^{m_i}, \]
where $\alpha$ runs over the same domain, where $v$ is regarded a column
vector, and where $\la .,. \ra$ denotes the natural pairing between column
vectors and row vectors. Thus $l_{i,b}$ is a piecewise linear
function, whose slopes correspond to the monomials in $f_{i,b}$.

\begin{thm} \label{thm:Join}
The dimension of $C_1+\ldots+C_k$ is at least the (topological) dimension
of the polyhedral set
\[ Q:=\left\{\left(\min_{i=1,\ldots,k} l_{i,b}(v_i)\right)_{b=1,\ldots,n} \mid 
	v_i \in \RR^{m_i} \text{ for all }
	i=1,\ldots,k\right\} \]
in $\RR^n$.
\end{thm}

As it stands, this theorem may not sound very appealing. For a more
concrete reformulation we proceed as follows. For $v=(v_1,\ldots,v_k)
\in \prod_{i=1}^k \RR^{m_i}$, $b \in \{1,\ldots,n\}$, and $i \in
\{1,\ldots,k\}$ we say that {\em $i$ wins $e_b$} (or $b$) {\em at $v$}
provided that
\begin{enumerate}
\item $l_{i,b}(v_i) < l_{j,b}(v_j)$ for all $j \neq i$, and 
\item $l_{i,b}$ is differentiable (hence linear) near $v_i$.
\end{enumerate}
If this is the case, then we denote by $d_{v_i} l_{i,b}$ the differential
$\RR^{m_i} \rightarrow \RR$ of $l_{i,b}$ at $v_i$; note that this is
given by a row vector of natural numbers. If, on the other hand, no $i$
wins $e_b$ at $v$, then we say that {\em there is a tie on $e_b$ at $v$}.

\begin{de}
For $v=(v_1,\ldots,v_k) \in \prod_{i=1}^k \RR^{m_i}$ and $i=1,\ldots,k$
set
\[W_i(v):=\{b \in \{1,\ldots,n\} \mid \text{$i$ wins $b$ at $v$}\}, \]
and call $W_i(v)$ the {\em winning set} of $i$ at $v$. Collect
the corresponding differentials $\RR^{m_i} \rightarrow \RR$ in the set
\[D_i(v):=\{d_{v_i} l_{i,b} \mid b \in W_i(v)\},\]
called the set of {\em winning directions} of $i$ at $v$. 
\end{de}

As we shall see in Section \ref{sec:Proofs}, the dimension of $Q$ in the
theorem is equal to the maximum, over all $v$, of
\[ \sum_{i=1}^k \dim_\RR \langle D_i(v) \rangle_\RR. \]
This leads to the following corollary.

\begin{cor} \label{cor:Join}
The dimension of $C_1+\ldots+C_k$ is at least the maximum, taken over
all $v=(v_1,\ldots,v_k) \in \prod_{i=1}^k \RR^{m_i}$, of the sum
\[ \sum_{i=1}^k \dim_\RR \langle D_i(v) \rangle_\RR. \]
\nopagebreak
In particular, if there exists a $v$ such that the set of winning
directions at $v$ of each $i$ spans a space of dimension $\dim C_i$
in $(\RR^{m_i})^*$, then the join $C_1+\ldots+C_k$ is non-defective.
\end{cor}

This corollary suggests the following strategy for proving that
$C_1+\ldots+C_k$ is non-defective: try and find a point $v$ at which
each $i$ wins a fair share of the basis $e_1,\ldots,e_n$---where fair
means that the linear forms on $\RR^{m_i}$ by means of which $i$ wins its
share, span a space of dimension $\dim C_i$. In the following section
we make this strategy more concrete for the case of secant varieties,
by making explicit the optimisation problem that needs to be solved to
get a good lower bound on $\dim kC$.

\section{Secant dimensions and some optimisation problems}
\label{sec:Secant}

Suppose that we want to find lower bounds on the secant dimensions of
a single closed cone $C \subseteq K^n$, which as before is the closure
of the image of a polynomial map $f=(f_1,\ldots,f_n): K^m \rightarrow
K^n=V$. For $b=1,\ldots,n$ let $A_b \subseteq \NN^m$ be the set of all
$\alpha$ for which the monomial $x^\alpha$ has a non-zero coefficient
in $f_b$. For $v=(v_1,\ldots,v_k) \in (\RR^m)^k$, $i=1,\ldots,k$, and
$b=1,\ldots,n$ we see that $i$ wins $e_b$ at $v$ if and only if there
is an $\alpha \in A_b$ which has $\la v_i,\alpha \ra < \la v_j,\beta
\ra$ for all $(j,\beta) \in \{1,\ldots,k\} \times A_b$ unequal to
$(i,\alpha)$. The winning set $W_i(v)$ of $i$ at $v$ is the set of all
$b$ with this property, and the set $D_i(v)$ of winning directions is
the set of all such minimising $\alpha$ as $b$ runs over $W_i(v)$. Hence
we are led to consider the following optimisation problem, in which we
relax the, at this point somewhat unnatural, restriction that all $A_b$
lie in $\NN^m$.

\begin{prb}[$\LP(A,k)$]
Let $A=(A_1,\ldots,A_n)$ be a sequence of finite subsets of $(\RR^m)^*$
and let $k \in \NN$. For $v=(v_1,\ldots,v_k) \in (\RR^m)^k$ and
$i=1,\ldots,k$ define
\begin{align*} D_i(v):=
	\bigcup_{b=1}^n 
	\{ \alpha \in A_b \mid 
		&\la v_i,\alpha \ra < \la v_j,\beta \ra \\
		&\text{ for all } (j,\beta) \in \{1,\ldots,m\} \times A_b
		\text{ unequal to } (i,\alpha) \}.
\end{align*}
{\bf Maximise} $\sum_{i=1}^k \dim \la D_i(v) \ra_\RR$ over all $v \in
(\RR^m)^k$.
\end{prb}

Note that at a given $v$, every $A_b$ only contributes to at most
one $D_i(v)$. The following corollary is immediate from Corollary
\ref{cor:Join}.

\begin{cor} \label{cor:Secant}
The dimension of $kC$ is at least the optimum of $\LP(A,k)$.
\end{cor}

The bad news is: this lower bound need not be very good. The good news: we
will see in Section \ref{sec:MinOrbit} that the bound is not useless---in
particular for minimal orbits in representations where all weight
spaces are one-dimensional---and in Section \ref{sec:Applications}
I give examples where the bound is very good. The following examples
illustrate both facts.

\begin{ex}[Bernd Sturmfels] \label{ex:Sturmfels}
Let $f:K^m \rightarrow K^n$ be a linear map whose matrix entries are
all non-zero. Then all $A_b$ are equal to $\{e_1^t,\ldots,e_m^t\}$,
where $e_j$ is the $j$-th standard basis vector of $\RR^m$.  Let
$v=(v_1,\ldots,v_k) \in (\RR^m)^k$. Now if $i$ wins $b$ at $v$, then
$i$ wins all other $b$ as well, and the winning directions $\alpha$
are all equal to the same $e_j^t$. Hence the lower bound of Corollary
$\ref{cor:Secant}$ on the dimension of $kC$ is $1$ for all $k \geq
1$---not very good indeed.
\end{ex}

\begin{ex} \label{ex:PolST}
Let $V=K[t_1,t_2]_d$ be the space of homogeneous polynomials of degree
$d$ in $t_1,t_2$, and let $C$ be the cone in $V$ of pure powers $(x_1
t_1+x_2 t_2)^d,\ x_1,x_2 \in K$. We prove that all secant varieties of
$C$ are non-defective \cite{Hirschowitz85}.  As a parametrisation
of $C$ we take the map $f:K^2 \rightarrow V,\ (x_1,x_2) \mapsto
(x_1t_1+x_2t_2)^d$. We fix the basis $t_1^b t_2^{d-b},\ b=0,\ldots,d$ of
$V$. The coefficient of $t_1^b t_2^{d-b}$ in $f(x_1,x_2)$ is $\binom{d}{b}
x_1^b x_2^{d-b}$, so that $A_b$ consists of the vector $(b,d-b)$
only. Hence Corollary \ref{cor:Secant} suggests that we compute the
optimum of $\LP(\{\{d-b\}\}_{b=0}^d, k)$.

\begin{figure}
\includegraphics[width=.5\textwidth]{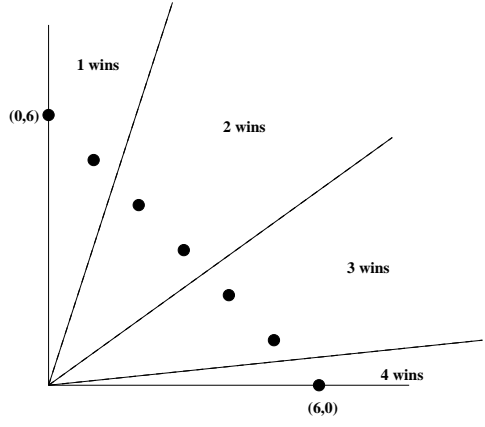}
\caption{The secants of $d$-th pure powers in binary forms are
non-defective.}
\label{fig:Veronese1}
\end{figure}

Now suppose first that $2k \leq d+1$. Then it is not hard to find a
$v=(v_1,\ldots,v_k)$ such that $D_{i+1}(v)$ contains both $(2i,d-2i)$
and $(2i+1,d-2i-1)$, for all $i=0,\ldots,k-1$: Draw $k-1$ lines
$l_1,\ldots,l_{k-1}$ in $\RR_{\geq 0}^2$, all through $0$ and such
that $l_i$ separates $(2i,d-2i)$ from $(2i-1,d-2i+1)$; and imagine any
concave piecewise linear (continuous) function $F$ on $\RR_{\geq 0}^2$
which is non-differentiable precisely in the $l_i$. Then the slope vectors
$v_i$ of $F$ on the $k$ components of the complement of the $l_i$ do the
job. So $D_{i+1}(v)$ contains two linearly independent vectors $(2i,d-2i),
(2i+1,d-2i-1)$, hence span has dimension $2=\dim C$. Using the corollary
we conclude that $kC$ is not defective.

If, on the other hand, $2k = d+2$, then one can find $v_1,\ldots,
v_{k-1}$ with the property above, while $k$ wins $(d,0)$; see Figure
\ref{fig:Veronese1}. Using the corollary we conclude that $kC=V$.
\end{ex}

Two variations of $\LP$ will appear in the sequel. First, there is
an affine version which is useful, for instance, when the map $f$ is
homogeneous, like in the preceding example. We use the notation $\Aff_\RR
D$ for the affine span of a subset $D$ in a real vector space. By
convention the dimension of $\emptyset=\Aff_\RR \emptyset$ is $-1$.

\begin{prb}[$\AP(A,k)$]
Let $A=(A_1,\ldots,A_n)$ be a sequence of finite subsets of $(\RR^m)^*$
and let $k \in \NN$. For $v=(v_1,\ldots,v_k) \in (\RR^m)^k$,
$a=(a_1,\ldots,a_k) \in \RR^k$, and $i=1,\ldots,k$, set
\begin{align*} 
D_i(v,a):= \bigcup_{b=1}^n 
	\{ \alpha \in A_b \mid 
	&\la v_i,\alpha \ra + a_i < \la v_j,\beta \ra +a_j\\
	&\text{ for all } (j,\beta) \in \{1,\ldots,m\} \times A_b
		\text{ unequal to } (i,\alpha) \}. 
\end{align*}
{\bf Maximise} $\sum_{i=1}^k (1+\dim \Aff_\RR D_i(v))$ over all $(v,a)
\in (\RR^m)^k \times \RR^k$.
\end{prb}

\begin{re} \label{re:AffineReduction}
The following obvious observation is sometimes useful: if
$A=(A_1,\ldots,A_n)$ is as in $\AP$, and $\pi:(\RR^m)^* \rightarrow
(\RR^{m'})^*$ is an affine-linear map, then the optimum of $\AP(A,k)$
is at least that of $\AP(\pi(A),k)$.
\end{re}

The third optimisation problem can be used when each component of the
parametrisation $f$ is a (non-zero multiple of a) single monomial,
like in Example \ref{ex:PolST}. This optimisation problem depends on
the choice of a positive definite inner product $(.,.)$ on $\RR^m$. We
use this inner product to identify $(\RR^m)^*$ with its dual $\RR^m$,
as well as to define a norm $||.||_2$ on $\RR^m$.

\begin{prb}[$\VP(S,k)$]
Let $S$ be a finite subset of $\RR^m$ and let $k \in \NN$. For
$v=(v_1,\ldots,v_k) \in (\RR^m)^k$ let $\Vor_i(v)$ denote the intersection
of $S$ with the Voronoi cell of $v_i$, i.e.,
\[
\Vor_i(v):=\{\alpha \in S \mid ||v_i-\alpha||_2 < ||v_j-\alpha||_2 
	\text{ for all } j \neq i \}.
\]
\noindent {\bf Maximise} $\sum_{i=1}^k (1+\dim \Aff_\RR \Vor_i(v))$ over
all $v \in (\RR^m)^k$.
\end{prb}

The relations between these optimisation problems are as follows.

\begin{lm} \label{lm:Interrelations}
Let $A=(A_1,\ldots,A_n)$ be a sequence of finite subsets of $(\RR^m)^*$,
and let $k \in \NN$.
\begin{enumerate}
\item If $\bigcup_b A_b$ is contained in an affine hyperplane not
through the origin, then the optimum of $\LP(A,k)$ equals the optimum
of $\AP(A,k)$.

\item If every $A_b$ is a singleton, then for $S:=\bigcup_b A_b$ the
optimum of $\VP(S,k)$ is a lower bound on the optimum of $\AP(A,k)$.

\end{enumerate}
\end{lm}

\begin{proof}
For the first statement: the affine-linear functions on $W:=\Aff_\RR
\bigcup_b A_b$ are precisely the restrictions to $W$ of the linear
functions on $\RR^m$; and furthermore $\dim \la M \ra_\RR=1+\dim \Aff_\RR
M$ for all $M \subseteq W$.

For the second statement, let $(v_1,\ldots,v_k)$ be an optimal solution
to $\VP(S,k)$; I will argue that there exist $v'=(v_1',\ldots,v_k') \in
(\RR^m)^k$ and $a=(a_1,\ldots,a_k) \in \RR^k$ such that
\begin{equation} \label{eq:VoronoiAffine}
||v_i-\alpha||_2 < ||v_j-\alpha||_2 \Rightarrow 
	(v_i',\alpha)+a_i < (v_j',\alpha)+a_j 
\end{equation}
for all $i,j=1,\ldots,m$ and $\alpha \in S=\bigcup_b A_b$.  This means
that the $\Vor_i(v)$ for $\VP(A,k)$ are contained in the $D_i(v',a)$
for $\AP(A,k)$, whence the lemma follows.

\begin{figure}
\includegraphics[width=.25\textwidth]{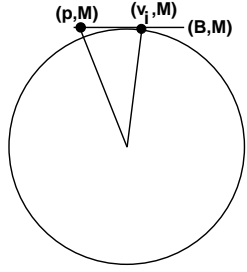}
\caption{The distance in $B \times \{M\}$ is close to the spherical
distance on a large sphere.}
\label{fig:Voronoi}
\end{figure}

Let $B$ be a compact convex set in $\RR^m$ containing the $v_i$ and
$S$; the latter is regarded as a subset of $\RR^m$ through the
inner product.  Embed $B$ in $\RR^{m+1}$ by giving every point in $B$ the
same $(m+1)$-st coordinate $M>0$; and extend the inner product to $\RR^{m+1}$
by making $e_{m+1}$ a norm-$1$ vector perpendicular to $\RR^m$. By making
$M$ large $B$ can be brought arbitrarily close to the sphere around the
origin in $\RR^{m+1}$ of radius $M$; see Figure \ref{fig:Voronoi}. In
particular, the function sending a point in $B$ to its distance to $v_i$
can be approximated, in the $\infty$-norm on continuous functions on $B$,
arbitrarily well by the spherical distance
\[ B \rightarrow \RR,\ 
x \mapsto M \arccos \frac{( (x,M), (v_i,M) )}{||(x,M)||_2 ||(v_i,M)||_2} \]
This, in turn, implies that the intersection with $B$ of the affine
hyperplane with equation $||x-v_i||_2=||x-v_j||_2$ can be arbitrarily well
approximated by the intersection with $B$ of the affine hyperplane
with equation
\[ \frac{((x,M),(v_i,M))}{||(v_i,M)||_2}
	= \frac{((x,M),(v_j,M))}{||(v_j,M)||_2}
	\left(= (x,v_j/||(v_j,M)||_2) + M^2/||(v_j,M)||_2 \right). \]
Hence, for $v'_i$ we take $-v_i / ||(v_i,M)||_2$, and for
$a_i$ we take $-M^2/||(v_i,M)||_2$; the minus signs
ensure that the $i$-th affine-linear function is the {\em minimal}
one near $v_i$, rather than the {\em maximal} one.  Then, for $M$
sufficiently large, \eqref{eq:VoronoiAffine} will be satisfied.
\end{proof}

Note that application of $\VP$ makes the proof in Example \ref{ex:PolST}
even easier: simply take $v_i$ in the middle between $(2i-2,d-2i+2)$
and $(2i-1,d-2i+1)$ for $i=1,\ldots,\lfloor \frac{d+1}{2} \rfloor$ and
$v_{1+\frac{d}{2}}$ equal to $(d,0)$ if $d$ is even; and note that this
$v$ gives the maximal possible value for $\VP$.

\section{Tropical geometry} \label{sec:Tropical}
Tropical geometry turns questions about algebraic varieties into questions
about polyhedral complexes, and this is precisely what the preceding
sections do to secant dimensions. For the general set-up, let
$K$ be an algebraically closed field of characteristic $0$, endowed with
a non-archimedean valuation $v:K \rightarrow \Rb:=\RR \cup \{\infty\}$,
which may, and in our application will, be trivial. Let $X$ be an affine
algebraic variety over $K$ and let $\lx=(\lx_1,\ldots,\lx_n)$ be an
$n$-tuple of generators of $K[X]$, giving rise to a closed embedding
$X \rightarrow K^n$. Let $x=(x_1,\ldots,x_n)$ be the standard coordinates
on $K^n$. In general, too, we will write $\bar{f} \in K[X]$ for the
restriction of a polynomial $f \in K[x]$ to $X$.

\begin{de}
The {\em tropicalisation of $X$ relative to $\lx$} is 
\[ T_\lx(X):=\left\{(v'(\lx_1),\ldots,v'(\lx_n)) \mid v':K[X]
	\rightarrow \Rb \text{ is a ring valuation extending
	}v\right\}. \] 
\end{de}

Here a {\em ring valuation extending $v$} is a map $v':K[X]
\rightarrow \Rb$ that equals $v$ on $K$ and satisfies the following
axioms: $v'(0)=\infty$, $v'(fg)=v'(f)+v'(g)$ and $v'(f+g) \geq
\min\{v'(f),v'(g)\}$ for all $f,g \in K[X]$. This definition, inspired
by \cite{Bieri84}, is the cleanest definition of the tropicalisation
of $X$. It shows clearly that $T_\lx(X)$ is only a projection of
an enormous object, namely the set of all ring valuations of $K[X]$
extending $v$. This fact explains why many things in tropical geometry
go almost right, but not quite; and why tropicalisation does not
have all the functorial properties one would like it to have. For
instance, if $X$ and $Y$ are both closed subvarieties in $K^n$, and
if we restrict the standard coordinates $x$ on $K^n$ to $X$, $Y$,
and $X\cap Y$, respectively, retaining the notation $\lx$ in all three
cases, then we obtain three tropicalisations $T_\lx(X)$, $T_\lx(Y)$,
$T_\lx(X \cap Y)$ in $\Rb^n$ and it is natural to ask: is $T_\lx(X
\cap Y)=T_\lx(X) \cap T_\lx(Y)$? While the inclusion $\subseteq$ is
clear from the definition, the converse inclusion does not always hold.
Another problem concerns morphisms: if $f:X \rightarrow Y$ is a morphism,
and if $\ly$ is a tuple of generators of $K[Y]$, then one would like to
have a tropicalisation $T(f): T_\lx(X) \rightarrow T_\ly(Y)$. Although
there is a natural such map if all pullbacks $f^*(\ly_j)$ 
are in the monoid generated by the $\lx_i$ (and in particular
$T(f)$ {\em can} be meaningfully defined on the tropical variety
$T_{\lx,f^*\ly}(X)$) in general $T(f)$ cannot be defined on all of
$T_\lx(X)$ in a meaningful way. For first steps in the abstract
theory of tropical varieties, see \cite{Mikhalkin06}.

Turning to applications of tropical geometry, we need a more useful
characterisation of $T_\lx(X)$. What follows stays very close to the
exposition in \cite{Speyer04}. I include it anyway for two reasons:
First, for self-containedness, and second, because there are some slight
differences: here we explicitly allow coordinates in $X$ to become $0$,
hence to have valuation $\infty$; and we make a clear distinction between
$K$ and a rather large valued extension $L$ that will soon play a role,
thereby emphasising that tropical methods even apply when the original
ground field is not endowed with a valuation.

For a $w \in \Rb^n, c \in K,$ and $\alpha \in \NN^n$ we call
$v(c)+\sum_i \alpha_i w_i$ the {\em $w$-weight} of the term $c x^\alpha
\in K[x]$, written $\wt_w c x^\alpha$; here we extend $+$ to $\Rb$ by
$a+\infty=\infty$ for all $a \in \Rb$ and set $0 \cdot \infty:=0$. For a
polynomial $f \in K[x]$ we let $\wt_w f$ be the minimum of the weights of
terms of $f$; in particular, $\wt_w 0=\infty$. Define the {\em $w$-initial
part} of $f$ to be $0$ if $\wt_w f=\infty$, and equal to the sum of all
terms in $f$ of weight $\wt_w f$ if the latter weight is $<\infty$.

Let $(L,v)$ be an algebraically closed and complete extension of $(K,v)$
with $v(L)=\Rb$ (such an extension exists.)  The definitions of $\wt_w$
and $\ini_w$ extend naturally to $L$.  The following theorem, which lies
at the heart of tropical geometry, gives an alternative description of the
tropicalisation $T_\lx(X)$. Both \cite{Einsiedler04} and \cite{Speyer04}
contain a version of this theorem; the proof below is close to that in
the latter reference---except that, like in the first reference, some
theory of affinoid algebras is used.

\begin{thm} \label{thm:Characterisations}
Let $I$ be the ideal of $X$ in $K[x]$. The following four sets are equal:
\begin{enumerate}
\item $\{ (v(\lx_1(p)),\ldots,v(\lx_n(p))) \mid p \in X(L) \}$,
\label{it:1}
\item $T_\lx(X)$, \label{it:2}
\item $\{ w \in \Rb^n \mid \ini_w f \text{ is not a monomial
	for any } f \in I\}$, \label{it:3} and 
\item $\{ w \in \Rb^n \mid \ini_w f \text{ is not a monomial
	for any } f \in L \otimes_K I\}$. \label{it:4}
\end{enumerate}
\end{thm}

This is really remarkable: the set of all ring valuations of $K[X]$
extending $v$ is huge and highly complicated. But this theorem says
that when one is only interested in the values of those valuations on
a {\em finite} tuple $\lx$ of generators of $K[X]$, then one needs only
consider the natural point valuations of $K[X]$ at $L$-rational points
of $X$. The proof uses the following two lemmas, both of which need
only slightly weaker assumptions on $L$ or $K$.

\begin{lm} \label{lm:Approximation}
Let $L$ be a field with a non-Archimedean valuation $v$ and let $K$ be a
subfield of $L$. Let $A$ be an $r \times s$-matrix with entries in $K$,
let $b \in K^r$, and let $l_1,\ldots,l_s$ be real numbers. Suppose that
there exists a $y \in L^s$ for which
\[ v((Ay-b)_i) > l_i \text{ for all } i=1,\ldots,r. \]
Then there also exists an $z \in K^s$ for which 
\[ v((Az-b)_i) > l_i \text{ for all } i=1,\ldots,r. \]
\end{lm}

First year linear algebra students know that if an {\em exact} solution
to the system $Az=b$ exists over $L$, then also one exists over $K$. This
lemma states that the same is true for {\em approximate} solutions. It
is not hard to prove the lemma using some theory of tensor products
of normed vector spaces, as contained in \cite[Chapter 1]{Fresnel04},
but here is an elementary proof.

\begin{proof}
As the statement only concerns the range of $A$, we may assume that
$A:K^s \rightarrow K^r$ is injective. In particular, we have $r \geq
s$, and we prove the lemma by induction on $r$. For $r=s$ the matrix
$A$ is invertible, so even an exact solution to $Az=b$ exists over
$K$. Now suppose that the statement is true for $r-1$, which is at least
$s$. Denote the rows of $A$ by $a_1,\ldots,a_r \in (K^s)^*$. As $r>s$,
there exists a linear relation $\sum_i \lambda_i a_i =0$ where not all
$\lambda_i$ are $0$. The existence of $y$ in the lemma yields
\begin{align}
v\left(\sum_i \lambda_i b_i\right)&= v\left(\sum_i \lambda_i (b_i-a_i y) + 
	\sum_i \lambda_i a_i y\right) \notag\\
	&=v\left(\sum_i \lambda_i (b_i-a_i y) + 0\right) \label{eq:lambdab}\\
	&>\min_i (v(\lambda_i)+l_i). \notag
\end{align}
After rearranging the rows of $A$ we may assume that the latter minimum is
attained in $i=r$, and by multiplying all $\lambda_i$ with $1/\lambda_r$ we
may assume that $\lambda_r=1$.	By the induction hypothesis, there exists
a $z \in K^s$ such that $v(a_i z -b_i)>l_i$ for all $i=1,\ldots,r-1$. 
For this same $z$ we have 
\begin{align*}
v(a_r z-b_r)&=v\left(-\sum_{i=1}^{r-1} \lambda_i a_i z -b_r\right)\\
&=v\left(-\sum_{i=1}^{r-1} \lambda_i (a_i z - b_i) - 
	\sum_{i=1}^r \lambda_i b_i \right)\\
&\geq \min \left\{v\left(\sum_{i=1}^{r-1} \lambda_i (a_i z - b_i)\right),
	v\left(\sum_{i=1}^r \lambda_i b_i \right)\right\}\\
&> \min_{i=1,\ldots,r} (v(\lambda_i)+l_i),
\end{align*}
where the last inequality follows from \eqref{eq:lambdab}
and the assumption on $z$. By assumption, the last minimum
is attained in $i=r$, and equal to $v(1)+l_r=l_r$.
\end{proof}

\begin{lm} \label{lm:Lifting}
Let $L$ be an algebraically closed field which is complete with
respect to a non-Archimedean valuation $v:L \rightarrow \Rb$. Set
$L^0:=\{c \in L \mid v(c) \geq 0\}$, $L^+:=\{c \in L \mid v(c)
>0\}$, and $\tilde{L}:=L^0/L^+$; the natural map $L^0 \rightarrow
\tilde{L}$, as well as all naturally induced maps, are denoted
$\pi$. Let $I$ be an ideal in $L[x_1,\ldots,x_n]$, set $I^0:=I \cap
L^0[x_1,\ldots,x_n]$ and $\tilde{I}:=\pi I^0$; the latter is an ideal
in $\tilde{L}[x_1,\ldots,x_n]$.

Then for any zero $\tilde{q} \in \tilde{L}^n$ of $\tilde{I}$ there exists
a zero $q \in (L^0)^n$ of $I$ for which $\pi(q)=\tilde{q}$.
\end{lm}

\begin{proof}
We prove this through an excursion to affinoid algebras; all properties
of these algebras that are used but not proved here can be found in
\cite{Bosch84, Fresnel04}. The motivation for such an excursion is the
following: let $T_n$ be the Tate algebra of all power series $\sum_{\alpha
\in \NN^n} c_\alpha x^\alpha$ over $L$ for which $\lim_{(\sum_i \alpha_i)
\rightarrow \infty} c_\alpha =0$. View $L[x_1,\ldots,x_n]$ as a subalgebra
of $T_n$ and let $J$ be the ideal in $T_n$ generated by $I$.  Then zeroes
of $I$ (and of $J$) in the polydisk $(L^0)^n$ correspond bijectively to
maximal ideals of the affinoid algebra $A:=T_n/J$.

To find such a zero we will use a lifting theorem from the theory of
affinoid algebras, whose formulation needs some further notions. First,
$T_n$ is a Banach algebra with the Gauss norm corresponding to the ring
valuation $v(\sum_\alpha c_\alpha x^\alpha)=\min_\alpha v(c_\alpha)$. As
any ideal of $T_n$, $J$ is closed---indeed, it is the closure of $I$---and
the affinoid algebra $A=T_n/J$ can be given the quotient norm, turning it
into a Banach algebra, as well. (In fact, any Banach algebra structure
on $A$ is equivalent to this one.) Let $A^0$ denote the subring of $A$
consisting of all $a$ for which the sequence $(a^n)_n$ is bounded, and
let $A^+$ denote the ideal in $A^0$ consisting of all $a$ with $\lim_{n
\rightarrow \infty} a^n=0$. One can show that $A^0$ maps any
zero $q$ of $J$ in $(L^0)^n$ into $L^0$, while $A^+$ maps it into
$L^+$, so that $q$ induces a $\tilde{L}$-algebra homomorphism
$\tilde{A}:=A^0/A^+ \rightarrow \tilde{L}$. We thus get a map
from $\{$zeroes of $J$ in $(L^0)^n\}$ to $\{\tilde{L}$-algebra
homomorphisms $\tilde{A} \rightarrow \tilde{L} \}$.  Theorem 3.5.3(ii)
in \cite{Fresnel04} says that this map is surjective; this is the lifting
theorem alluded to before.

To use this theorem, we must construct a homomorphism $\tilde{A}
\rightarrow \tilde{L}$ from our point $\tilde{q}$. This goes as follows:
first let $q'\in (L^0)^n$ be {\em any} lift of $\tilde{q}$. By assumption
$I^0$ maps $q'$ into $L^+$ and hence, since $J^0:=\{j \in J \mid v(j) \geq
0\}$ is the closure of $I^0$ and $L^+$ is closed, $J^0$ also maps $q'$
into $L^+$. We conclude that the map $\phi: T^0/J^0 \rightarrow \tilde{L}$
sending $f+J^0$ to $\pi(f(q'))$ is a well-defined ring homomorphism,
which restricts to $\pi$ on $L^0$.  Since $A^0$ is integral over its
subring $T^0/J^0$ \cite[Theorem 3.5.3(i)(2)]{Fresnel04} we can find
a prime $Q$ in $A^0$ lying over the prime $\ker \phi$, so that the
following diagram commutes:
\[ 
\xymatrix{
T^0/J^0 \ar[r] \ar@{->>}[d]_\phi \ar@{^{(}->}[r] & A^0 \ar@{->>}[d]\\
\tilde{L} \ar@{^{(}->}[r] & A^0/Q
}
\]
Now the lower horizontal inclusion is again integral, and since $\tilde{L}$
is an algebraically closed field it is an isomorphism. In other
words, $\phi$ factorises through a ring homomorphism $A^0 \rightarrow
\tilde{L}$---which on $L^0$ still equals $\pi$, of course. From $A^+=L^+
A^0$ and $\phi(L^+)=0$ we conclude that $A^+$ is mapped to zero under this
homomorphism, and we have the desired $\tilde{L}$-algebra homomorphism
$\psi: \tilde{A} \rightarrow \tilde{L}$.

By the aforementioned theorem, we can now find a zero $q \in (L^0)^n$
of $J$ such that $\pi(f(q))=\psi(f+A^+)$ for all $f \in A^0$. Taking
for $f$ the coordinate function $x_i + J_0 \in T^0/J^0 \subseteq A^0$ we
find that $\pi(q_i)=\psi(x_i+J_0+A^+)=\phi(x_i+J_0)=\tilde{q}_i$, as desired.
\end{proof}

\begin{proof}[Proof of Theorem \ref{thm:Characterisations}]
The inclusion \eqref{it:1}$\subseteq$\eqref{it:2} follows from the
fact that, for a point $p \in X(L)$, the map $K[X] \rightarrow \Rb,\
\bar{f} \mapsto v(\bar{f}(p))$ is a ring valuation extending $v$. For
the inclusion \eqref{it:2}$\subseteq$\eqref{it:3} let $w \in T_\lx(X)$
and let $v':K[X] \rightarrow \Rb$ be a ring valuation extending $v$
with $v'(\lx_i)=w_i$.  Let $f=\sum_\alpha c_\alpha x^\alpha \in I$. We
show that $\ini_w f$ is not a monomial: Indeed, suppose that
it is. Then on the one hand, $v'(\bar{f})=v'(0)=\infty$ by the first axiom
for ring valuations, while on the other hand $v'(\bar{f})= v'(\sum_\alpha
c_\alpha \lx^\alpha)$. Now the value $v'(c_\alpha \lx^\alpha)$ equals
$\wt_w c_\alpha x^\alpha$ by definition of $\wt$ and the choice of
$v'$, and by assumption a unique term of $f$ minimises this value.
The axioms of ring valuations readily imply that $v'(\bar{g}+\bar{h})
= v'(\bar{g})$ whenever $v'(\bar{g})<v'(\bar{h})$, hence $v'(\bar{f})$
equals that uniquely assumed minimal value, which is smaller than
$\infty$---a contradiction to $v'(\bar{f})=\infty$. 

For the inclusion \eqref{it:3}$\subseteq$\eqref{it:4}, let $w \in \Rb^n$
be such that $\ini_w f$ is not a monomial for any $f \in I$ and suppose
that $\ini_w \sum_{j=1}^s c_j f_j$ {\em is} a monomial $g$ for some $c_j
\in L$ and $f_j \in I$.  Let $g_1:=g,g_2,\ldots,g_r$ be the collection of all
monomials occurring in the $f_i$, and set $l_i:=\wt_w(g)-\wt_w(g_i)$
for $i=1,\ldots,r$. Let $A$ be the $r \times s$-matrix over $K$ whose
entry at position $(i,j)$ is the coefficient of $g_i$ in $f_j$. Then
the vector $c:=(c_1,\ldots,c_s)^t \in L^s$ satisfies $v((Ac-e_1)_i)>l_i$
for all $i=1,\ldots,r$; hence by Lemma \ref{lm:Approximation} there also
exists a $c':=(c_1',\ldots,c_s')^t \in K^s$ for which $v((Ac'-e_1)_i)>l_i$
for all $i$. But then $\ini_w \sum_i c_i' f_i \in I$ is a non-zero scalar
multiple of $g=g_1$, as well---a contradiction to \eqref{it:3}.

Finally, for \eqref{it:4}$\subseteq$\eqref{it:1}, let $w \in \Rb^n$ be
such that $\ini_w f$ is not a monomial for any $f \in L \otimes_K I$.
We have to exhibit a point $p=(p_1,\ldots,p_n) \in X(L)\subseteq
L^n$ with $v(p_i)=w_i$ for all $i$. Choose $t_1,\ldots,t_n \in L$
with $v(t_i)=w_i$ for all $i$, set $x':=(x_i)_{i:w_i \neq \infty}$,
and let $\phi:L[x] \rightarrow L[x']$ be the epimorphism sending $f$
to $f(t_1x_1,\ldots,t_nx_n)$; $\phi(f)$ only contains variables $x_i$
with $w_i \neq \infty$ since the remaining ones have $t_i=0$. We then
have $\wt_w f=\wt_0 \phi(f)$ and $\phi(\ini_w f)=\ini_0 \phi(f)$. Let
$J\subseteq L[x']$ be the ideal $\phi(L \otimes_K I)$; by construction
$\ini_0 f$ is not a monomial for any $f \in J$. We claim that there exists
a $q \in L^{\{i:w_i \neq \infty\}}$ such that $v(q_i)=0$ for all $i$
and such that $J$ vanishes on $q$; then setting $p_i:=t_iq_i$ if $w_i
\neq \infty$ and $p_i:=0$ if $w_i=\infty$ gives a point $p$ as required.

Now retain the notation $L^0,L^+,\tilde{L},\pi$ from Lemma
\ref{lm:Lifting}. Let $J^0:=J \cap L^0[x']$, and and set $\tilde{J}:=\pi
J^0$. Since for any $f \in J^0$ either $\pi(f)$ is $0$ or $\pi(f)$
has the same monomials as $\ini_0 f$, the ideal $\tilde{J} \subseteq
\tilde{L}[x']$ contains no monomials. As $\tilde{L}$ is algebraically
closed, $\tilde{J}$ has a zero $\tilde{q}$ in $(\tilde{L}^*)^{\{i:w_i
\neq \infty\}}$ by the Nullstellensatz. Applying Lemma \ref{lm:Lifting}
to $J \subseteq L[x']$, we conclude that $\tilde{q}$ can be lifted to a
zero $q \in (L^0)^{\{i:w_i \neq \infty\}}$ of $J$. Clearly all components
of $q$ have valuation $0$, so we are done.

\end{proof}

By Theorem \ref{thm:Characterisations}, $T_\lx(X)$ is the intersection
of infinitely many polyhedral sets, one for each element $f$ of $I$:
the set of all $w \in \Rb^n$ for which $\ini_w f$ is not a monomial. One
can show that, in fact, finitely many of these polyhedral sets already
cut out $T_\lx(X)$ \cite{Sturmfels02}, so that the latter set is a
polyhedral complex.  The following theorem, originally due to Bieri
and Groves \cite{Bieri84} and also proved in \cite{Sturmfels02} using
Gr\"obner basis methods, relates the dimension of this polyhedral set
to that of $X$.

\begin{thm} \label{thm:BieriGroves}
Suppose that $X$ is irreducible and of dimension $d$. Then $T_\lx(X)$ is a
polyhedral complex in $\Rb^n$ which is pure of dimension $d$.
\end{thm}

As mentioned before, there is no obvious tropicalisation of morphisms
between embedded affine varieties. However, {\em polynomial
maps} do have natural tropicalisations. 

\begin{de}
For a polynomial $h \in K[x_1,\ldots,x_m]$, the map
\[ T(h):\Rb^m \rightarrow \Rb,\ w \mapsto \wt_w h \]
is called the {\em tropicalisation} of $h$. Similarly, for a
polynomial map $f=(f_1,\ldots,f_n):K^m \rightarrow K^n$, the
map
\[ T(f):\Rb^m \rightarrow \Rb^n,\ T(f):=(T(f_1),\ldots,T(f_n)) \]
is called the {\em tropicalisation} of $f$.
\end{de}

Note that $T(f)$ is continuous when we give $\Rb$ the usual topology of
a half-open interval. The following lemma is also well-known; see for
instance \cite[Theorem 3.42]{Pachter05} for a more detailed statement.
I include its short proof for self-containedness.

\begin{lm} \label{lm:TropMap}
Let $f:K^m \rightarrow K^n$ be a polynomial map, let $X$ be the
Zariski closure of $\im(f)$, and let $\lx=(\lx_1,\ldots,\lx_n)$ be the
restrictions to $X$ of the standard coordinates on $K^n$. Then $T(f)$
maps $\Rb^m$ into $T_\lx (X)$.
\end{lm}

\begin{proof}
First let $w \in \RR^n$ be such that for all $i=1,\ldots,n$ either $f_i$
is identically zero or $\ini_w f_i$ is a monomial; note that then $T(f)$
is affine-linear near $w$ (disregarding the infinite entries corresponding
to the $f_i$ that are identically $0$). Choose $p=(p_1,\ldots,p_m) \in L^m$
with $v(p_i)=w_i$ for all $i$. Then the definition of $T(f_i)$ 
implies that $v(f_i(p))= T(f_i)(w)$ for all $i$. Hence $T(f)$ maps
such a $w$ into $T_\lx(X)$ by the inclusion \eqref{it:1} $\subseteq$
\eqref{it:2} of Theorem \ref{thm:Characterisations}. The set of all
$w$ with the required properties is clearly dense in $\Rb^n$, $T(f)$
is continuous, and $T_\lx(X)$ is closed---whence the lemma.
\end{proof}

\section{Proofs of Theorem \ref{thm:Join} and Corollary
	\ref{cor:Join}} \label{sec:Proofs}

We retain the notation of Section \ref{sec:Join}.

\begin{proof}[Proof of Theorem \ref{thm:Join}.]
Consider the map $f:\prod_{i=1}^k K^{m_i} \rightarrow C_1+\ldots+C_k
\subseteq V$ sending $(p_1,\ldots,p_k)$ to $f_1(p_1)+\ldots+f_k(p_k)$.
Endow $K$ with the trivial valuation, and observe that the map
$\prod_{i=1}^k \RR^{m_i} \rightarrow \RR^{m}$ whose $b$-th component
is $\min_{i=1,\ldots,k} l_{i,b}$ is precisely the tropicalisation
$T(f)$ of $f$---or rather, its restriction to the set $\prod_{i=1}^k
\RR^{m_i}$. Hence the set $Q$ in the theorem is precisely $T(f)(\prod_{i=1}^k
\RR^{m_i})$. By Lemma \ref{lm:TropMap} this set is contained in the
tropicalisation $T_\ly (C_1+\ldots+C_k)$ (where $y=(y_1,\ldots,y_n)$
are the standard coordinates on $K^n$), and hence by Theorem
\ref{thm:BieriGroves} its dimension does not exceed the dimension of
$C_1+\ldots+C_k$.
\end{proof}

\begin{proof}[Proof of Corollary \ref{cor:Join}.]
Let $f$ be the polynomial map from the previous proof, so that $T(f)$ is
a piecewise linear map $\prod_{i=1}^k \RR^{m_i} \rightarrow \RR^n$. Let
$v \in \prod_i \RR^{m_i}$ and let $B \subseteq \{1,\ldots,n\}$ be the
set of indices won at $v$ by some $i \in \{1,\ldots,k\}$; that is, we
simply leave out the indices where there is a tie. Then the map $F$,
defined as the composition of $T(f):\prod_i \RR^{m_i} \rightarrow \RR^n$
and the projection $\RR^n \rightarrow \RR^B$, is linear near $v$,
so that its differential $d_v F:\prod_i \RR^{m_i} \rightarrow \RR^B$
is well-defined. Indeed, for $i=1,\ldots,k$ and $b \in B$ the $b$-th
component of the restriction of $d_v F$ to $\RR^{m_i}$ is $d_{v_i}
l_{i,b}$ if $i$ wins $b$, and $0$ otherwise. Hence the rank of $d_v T(f)$
is exactly
\[ \sum_{i=1}^k \dim_\RR \langle D_i(v) \rangle_\RR; \] 
therefore this number is a lower bound on $\dim \im F$, which in turn
is a lower bound on $\dim \im T(f)$, and hence to $\dim(C_1+\ldots+C_k)$
by the proof of Theorem \ref{thm:Join}.
\end{proof}

Note that the dimension of $Q$ in Theorem \ref{thm:Join} is in fact equal
to the maximum rand of $d_v T(f)$ over all a points $v$ where this rank
is linear.

\section{Minimal orbits} \label{sec:MinOrbit}
An important class of varieties, whose non-defectiveness is notoriously
hard to prove, are the {\em minimal orbits}. To define these, let $G$
be a connected reductive algebraic group over $K$, and let $V$ be an
irreducible $G$-module. Then $\PP V$ has a unique closed $G$-orbit
(see for instance \cite{Borel91, Humphreys72, Humphreys81} for general
theory of algebraic groups, Lie algebras, and their representations),
the cone over which is $C:=G v_\lambda \cup \{0\} \subseteq
V$, where $v_\lambda$
is the highest weight vector of $V$ relative to some Borel subgroup
$B$. Many interesting cones $C$ arise in this manner.

\begin{ex} \label{ex:Minorbits}
\begin{enumerate}
\item Let $G$ be $\lieg{GL}_m$, let $d \in
\{1,\ldots,m-1\}$, and set $V$
equal to $\Wedge^d (K^m)$. Then $C$ is the cone over the Grassmannian,
in its Pl\"ucker embedding, of $d$-dimensional subspaces of $K^m$.

\item Let $G$ be $(\lieg{GL}_m)^d$, let $d \in \NN$, and let $V$ be the space
$(K^m)^{\otimes d}$. Then $C$ is the cone of all {\em pure $d$-tensors},
i.e., those that can be written as $u_1 \otimes \ldots \otimes u_d$
for some $u_1,\ldots,u_d \in K^m$; it is the cone over the
Segre embedding of the $d$-fold Cartesian power of $\PP^{m-1}$.

\item Let $G$ be $\lieg{GL}_m$, let $d \in \NN$, and let $V$ be the
space of homogeneous polynomials of degree $d$ on $K^m$. Now $C$
is the cone of all pure $d$-th powers of linear forms, and the affine
cone over the Veronese embedding of $\PP^{m-1}$ of degree $d$. It is
well known that the dimension of $kC$ is the codimension of the space
of all homogeneous polynomials $f$ of degree $d$ for which both $f$
and all first order partial derivatives of $f$ vanish in $k$ generic
points on $\PP^{m-1}$. This relates the secant dimensions of $C$ to
{\em multivariate polynomial interpolation}. These dimensions are known
for all $m$ and $d$ from the papers \cite{Alexander88, Alexander92,
Alexander95, Hirschowitz85}.

\item Let $G$ be $\lieg{SO}_m$, and let $V$ be the Lie algebra of
$G$. Then $C$ is the cone over the Grassmannian, in its Pl\"ucker
embedding, of isotropic $2$-dimensional subspaces of $K^m$. The secant
varieties of $C$ were determined in \cite{Baur04}.

\end{enumerate}
\end{ex}

There is a general argument showing that the first few secant varieties of
minimal orbits are non-defective. To state it, let $T$ be a maximal torus
in $B$, let $P \supseteq B$ be the stabiliser of $Kv_\lambda$, and let $U$ be
the unipotent radical of the parabolic subgroup opposite to $P$ containing $T$.
Then $\liea{g} v_\lambda=\liea{u} v_\lambda \oplus
Kv_\lambda$. Finally, let $W$ be the
Weyl group of $(G,T)$; for $w \in W$ let $\tilde{w}$ be an element of
$N_G(T)$ representing $w$.

\begin{prop} \label{prop:WeylOrbit}
Let $k \in \NN$ and $w_1,\ldots,w_k$ be elements of $W$. Then
\[ \dim kC \geq \dim \sum_{i=1}^k \tilde{w}_i
(\liea{u}v_\lambda +
Kv_\lambda). \]
In particular, if there exist $w_1,\ldots,w_k$ for which the spaces
$\tilde{w}_i(\liea{u}v_\lambda + Kv_\lambda)$ are linearly independent, then $kC$
is non-defective.
\end{prop}

\begin{proof}
The rank of the differential of the addition map $C^k \rightarrow
kC$ at any point of $C^k$ is a lower bound on the dimension of $kC$,
which in an open dense subset of $C^k$ is exact---this is Terracini's
lemma \cite{Terracini11}. Now take for this point the point $(\tilde{w}_1
v_\lambda,\ldots,\tilde{w}_kv_\lambda)$. The tangent space to $C$ at $\tilde{w}_i
v_\lambda$ is $\tilde{w}_i
\liea{g}v_\lambda=\tilde{w}_i(\liea{u}v_\lambda + Kv_\lambda)$, and
the differential of the summation map maps the tuple of these spaces to
their sum.
\end{proof}

Proposition \ref{prop:WeylOrbit} is useful for small $k$ and large highest
weights: the space $\liea{u}v_\lambda + Kv_\lambda$ is contained in the sum of the
weight spaces with weights $\lambda - \alpha$, where $\lambda$ is the
weight of $v_\lambda$ and $-\alpha$ is $0$ or a root whose root space lies in
$\liea{u}$. Hence if there exist $w_1,\ldots,w_k$ such that the translates
\[ w_i \{\lambda-\alpha \mid -\alpha \text{ is a root of }
	\liea{u} \} \]
are all disjoint, then $kC$ is non-defective. For $\lambda$ large and
deep in the interior of the dominant chamber, there will exist such
$w_i$ for all $k$ up to $|W|$. However, for $k>|W|$, the bound of the
proposition is evidently off.

The bound from Proposition \ref{prop:WeylOrbit} is also quite good for
representations that are small in the sense that all weights, or many
of them, are highest weights. Let me illustrate this in the first three
examples of \ref{ex:Minorbits}. For the coding theory notions appearing
here and in what follows, I refer to \cite{Lint99}.

\begin{ex} \label{ex:WeylOrbit}
\begin{description}
\item[Grassmannians] Take for $T$ the diagonal matrices.  The $W$-orbit
of the highest weight vector consists of the vectors $e_{i_1}\wedge
\ldots \wedge e_{i_d}$ with $1 \leq i_1<\ldots<i_d \leq m$. These vectors
correspond naturally to the binary words in $\{0,1\}^m$ of Hamming weight
$d$. Given a collection $B$ of $k$ such words, the lower bound of the
proposition equals the number of weight-$d$, length-$m$ binary words
at distance at most $2$ from $B$. Hence, if there exists a binary code
$B$ of size $k$, length $m$, constant weight $d$, and minimal Hamming
distance distance $6$, then $kC$ is non-defective.  Variants of this
idea already appeared in \cite{Ehrenborg00,Catalisano05}.

\item[Pure tensors] Take for $T$ the $d$-tuples of diagonal matrices. The
$W$-orbit of the highest weight vector consists of all tensors of the
form $e_{i_1} \otimes \ldots \otimes e_{i_d}$, where $i_1,\ldots,i_d
\in \{1,\ldots,m\}$. These correspond naturally to the $m$-ary words
of length $d$. Given a collection $B$ of such words, the lower bound
from the proposition equals the number of words at Hamming distance at
most $1$ from $B$, or, equivalently, to the total number of fields on
a $d$-dimensional chessboard with side lengths $m$ covered by rooks
on the positions in $B$. Hence, if there exists an $m$-ary code $B$
of size $k$, length $d$, and minimal Hamming distance $3$, then $kC$ is
non-defective. The existence of perfect (Hamming) codes for $m=p^e$
and $d=(p^{ef}-1)/(p^e-1)$, where $p$ is a prime and $e$ and $f$ are
arbitrary positive natural numbers, shows that the corresponding $C$
are non-defective. This idea is contained in \cite{Catalisano2002,
Ehrenborg00, Sturmfels05}.

\item[Homogeneous polynomials] Take for $T$ the diagonal matrices. Here
the weight vectors are the monomials, corresponding naturally to
the multi-indices $\alpha \in \NN^m$ with $|\alpha|:=\sum_{i=1}^m
\alpha_i=d$. The elements of the form $de_i$ correspond to the elements
$x_i^d$ in the $W$-orbit of the highest weight vector. If $B$ is a
set of such $de_i$, then the lower bound from the proposition is the
number of $\alpha$ that are at $1$-distance at most $2$ from $B$; this
is also mentioned in \cite{Ehrenborg00}.  Dually, if $B$ has size $k$
and consists of elements that are mutually at $1$-distance at least $6$,
then $kC$ is non-defective.  Note that such a $B$ exists if and only if
$k \leq m$ and $d \geq 3$.

For $d=2$ the bound turns out to be exact: Taking
$B=\{2e_1,\ldots,2e_k\}$, where $k \leq m$, the proposition says that
$kC$ has dimension at least $m+(m-1)+\ldots+(m-k+1)$, and this is known
to be the exact dimension.
\end{description}
\end{ex}

It would be very interesting to apply Corollary \ref{cor:Secant} to
minimal orbits in general irreducible representations, but for this
one needs a suitable basis of the representation $V$, and it is not
obvious which one to use. This will be the theme of further research;
here we content ourselves with the following lower bound, which seems
to be good for Segre products of Veronese embeddings. Write $X(V)$
for the set of $T$-weights on $V$, considered as a subset of the real
vector space $\RR \otimes_\ZZ X(T)$ spanned by the character group of $T$.

\begin{prop} \label{prop:MultOne}
For any positive definite inner product on $\RR \otimes_\ZZ X(T)$ and
any $k$, the optimum of $\VP(X(V),k)$ is a lower bound on $\dim kC$.
\end{prop}

\begin{proof}
To apply Corollary \ref{cor:Secant} we need a dominant polynomial map
into $C$ together with choices of basis. For the map we take
\[ f:K \times \liea{u} \rightarrow C \subseteq V,\ (t,u) \mapsto
	t\exp(u)v_\lambda. \]
Let $X(\liea{u})$ denote the set of $T$-roots in $\liea{u}$ and set
$\tilde{X}(\liea{u}):=\{0\} \cup X(\liea{u})$. In $K \times \liea{u}$ we
choose a basis labelled by $\tilde{X}(\liea{u})$, where $0$ corresponds
to $(1,0)$ and $\beta \in X(\liea{u})$ corresponds to a root vector
$u_\beta$ in $\liea{u}$ with root $\beta$. Let $(x_\beta)_{\beta \in
\tilde{X}(\liea{u})}$ be the corresponding coordinates on $K\times
\liea{u}$.  Finally, in $V$ we choose any basis $v_1,\ldots,v_n$ of
$T$-weight vectors.

For $b=1,\ldots,n$ let $A_b$ be the set of exponent vectors of
monomials in the $x_\beta,\beta \in \tilde{X}(\liea{u})$, occurring
in $f_b$. Observe that if $\mu$ is the weight of $v_b$, then these
exponent vectors $(r_\beta)_{\beta \in \tilde{X}(\liea{u})}$ all satisfy
$r_0=1$ and furthermore $\lambda+\sum_\beta r_\beta \beta=\mu$---only
such monomials in the $u_\beta$ can map $v_\lambda$ to an element
having non-zero $b$-th component. By the latter equality, Corollary
\ref{cor:Secant}, and Lemma \ref{lm:Interrelations}, $\AP((A_b)_b,k)$
is a lower bound on $\dim kC$.

Now let $\pi$ be the affine-linear map from
$(\RR^{\tilde{X}(\liea{u})})^*$ to $\RR \otimes_\ZZ X(T)$ sending
$(r_\beta)_\beta$ to $\lambda+\sum_\beta r_\beta \beta$; the above shows
that $\pi(A_b)=\{\mu\}$ if $v_b$ has weight $\mu \in X(V)$.  Hence,
by Remark \ref{re:AffineReduction} and Lemma \ref{lm:Interrelations},
$\AP((\{\mu\})_\mu,k)$ and $\VP(X(V),k)$ are lower bounds to $\dim kC$,
as well.
\end{proof}

Clearly, if $V_\mu$ is $1$-dimensional for all $\mu$, then the bound
of Proposition \ref{prop:MultOne} equals $\dim V$ for $k$ sufficiently
large---a property that did not hold for the bound of Proposition
\ref{prop:WeylOrbit}. This shows that the bound is not {\em useless},
as in Example \ref{ex:Sturmfels}. I do not claim that the bound
of Proposition \ref{prop:MultOne} is {\em good} for all weight-$1$
representations. For instance, for $C$ the cone over the Grassmannian of
$d$-dimensional subspaces of an $m$-dimensional subspace, the bound only
gives $\dim C \geq m$ rather than $d(m-d)+1$. This is due to the rather
coarse reduction using $\pi$; so to apply Corollary \ref{cor:Secant}
to Grassmannians one needs a more subtle approach. When $X(\liea{u})$
is linearly independent, however, as is the case for Veronese and Segre
embeddings and indeed for Segre products of Veronese embeddings, then
$\pi$ is an isomorphism, and the bound of Proposition \ref{prop:MultOne}
seems actually very good, as we will see in the next
section. 

There is much more to be said here: what about representations where
the weight multiplicities are not $1$?  What about representations
where they are, but $X(\liea{u})$ is linearly dependent? Such questions
will be addressed in forthcoming research; I conclude this paper with a
few concrete examples of how strong the lower bounds from the tropical
approach can be.

\section{Applications} \label{sec:Applications}

\subsection{Veronese embeddings}
Let $m,d$ be natural numbers and set $V:=\CC[t]_d$, the space of
complex homogeneous polynomials of degree $d$ in the variables
$t=(t_1,\ldots,t_m)$. As a basis of $V$ we choose the monomials
$t^\alpha$, where $\alpha$ runs through the set $X \subseteq \NN^m$
of all multi-indices with $|\alpha|=d$. Let $f:K^m \rightarrow V$
send $(x_1,\ldots,x_m)$ to the pure power $(x_1 t_1+\ldots+x_mt_m)^d$,
and let $C$ be the image of $f$. The coefficient of $t^\alpha$ in
$f(x)$ is a non-zero scalar times $x^\alpha$, so---like in Example
\ref{ex:PolST}, where $m$ was $2$---every $A_\alpha$ consists of the
single vector $\alpha$.  So Corollary \ref{cor:Secant} leads us to
consider the optimisation problem $\LP(\{\{\alpha\}\}_{\alpha \in X},k)$
to bound the dimension of $kC$. As $X$ lies on an affine hyperplane
not through $0$, we may by Lemma \ref{lm:Interrelations} just as well
solve $\AP(\{\{\alpha\}\}_{\alpha \in X},k)$, and since all sets $A_b$
are singletons, also $\VP(X,k)$ gives a lower bound on $\dim kC$; this is
also the content of Proposition \ref{prop:MultOne} for the representation
of $\lieg{GL}_m$ on $V$.

Of course we know the dimension of $kC$ already from \cite{Alexander95},
but it would be very desirable to have an alternative, more elementary
proof of their results. I think that the tropical approach might yield
such a proof. To motivate this belief, let us prove the first non-trivial
case solved by Hirschowitz in \cite{Hirschowitz85}, namely, the case
where $m=3$. Rick Miranda and Olivia Dumitrescu have also proved the
following theorem, using degenerations (private communication). The
combinatorics to which their proof boils down resembles very much the
combinatorics below, and it would be interesting to understand exactly
how both approaches are connected.

\begin{thm} \label{thm:Veronesem3}
Suppose that $m=3$. Then $kC$ has the expected dimension
unless $(d,k)=(2,2)$ or $(4,5)$, in which cases the defect 
$k\dim C-\dim kC$ equals $1$.
\end{thm}

\begin{figure}
\subfigure[$d=2$]{\includegraphics[width=.2\textwidth]{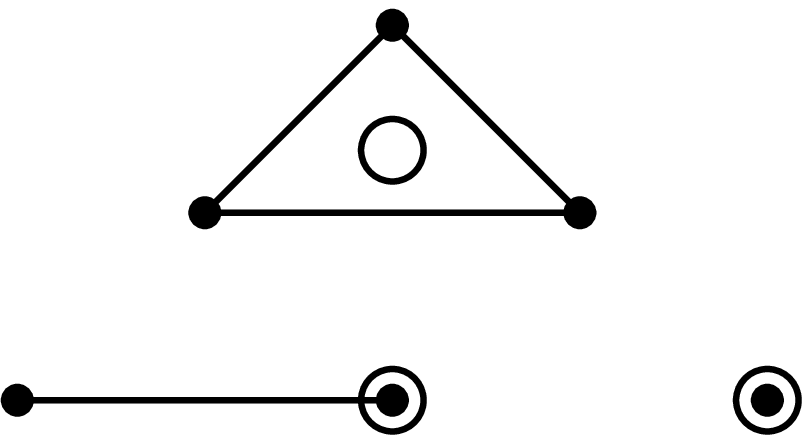}\label{fig:m3d2}}
\hfill
\subfigure[$d=3$]{\includegraphics[width=.3\textwidth]{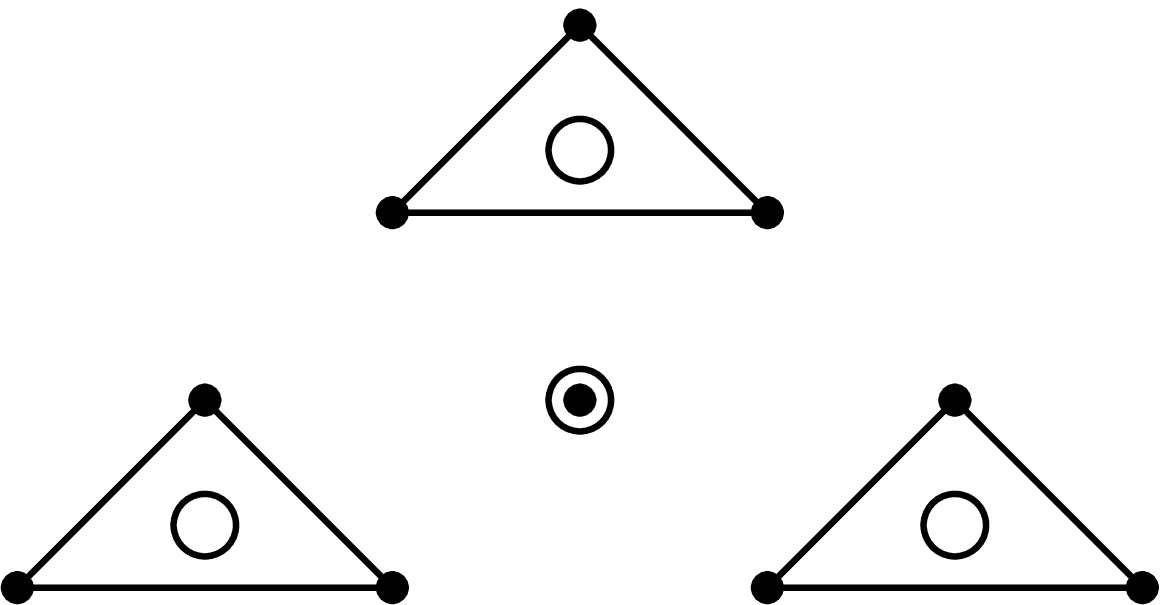}\label{fig:m3d3}}
\hfill
\subfigure[$d=4$]{\includegraphics[width=.3\textwidth]{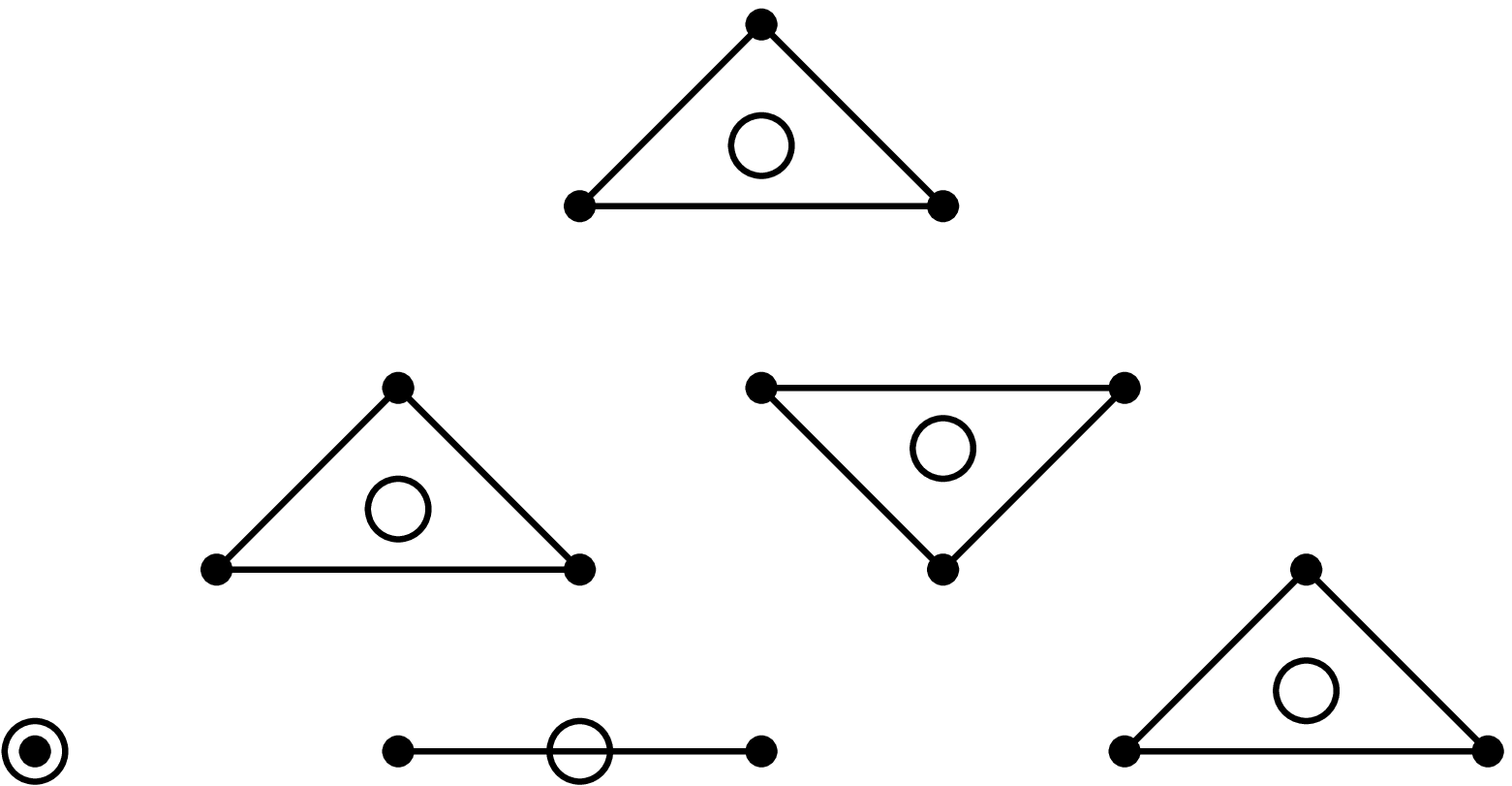}\label{fig:m3d4}}
\hfill
\subfigure[$d=5$]{\includegraphics[width=.4\textwidth]{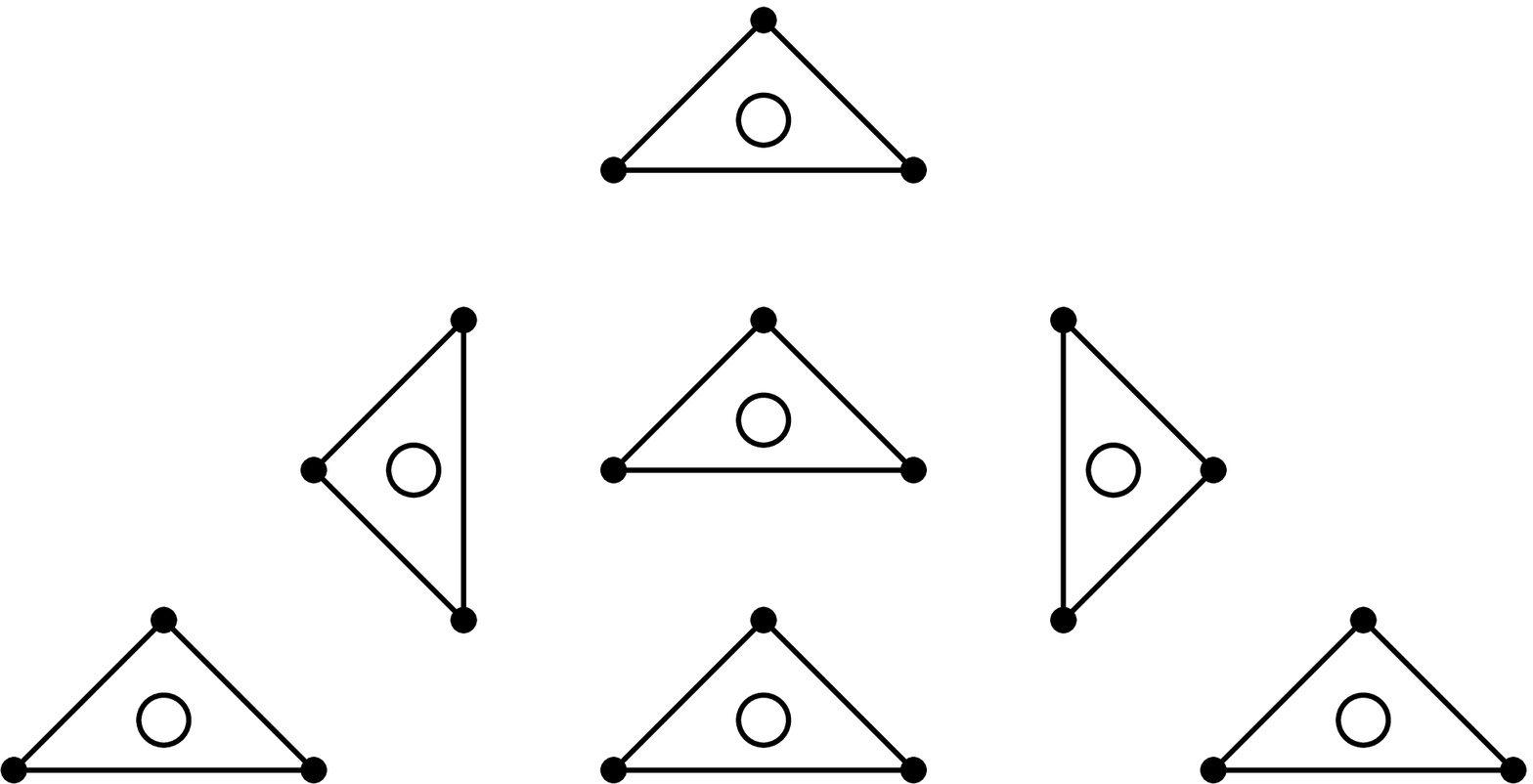}\label{fig:m3d5}}
\hfill
\subfigure[$d=6$]{\includegraphics[width=.4\textwidth]{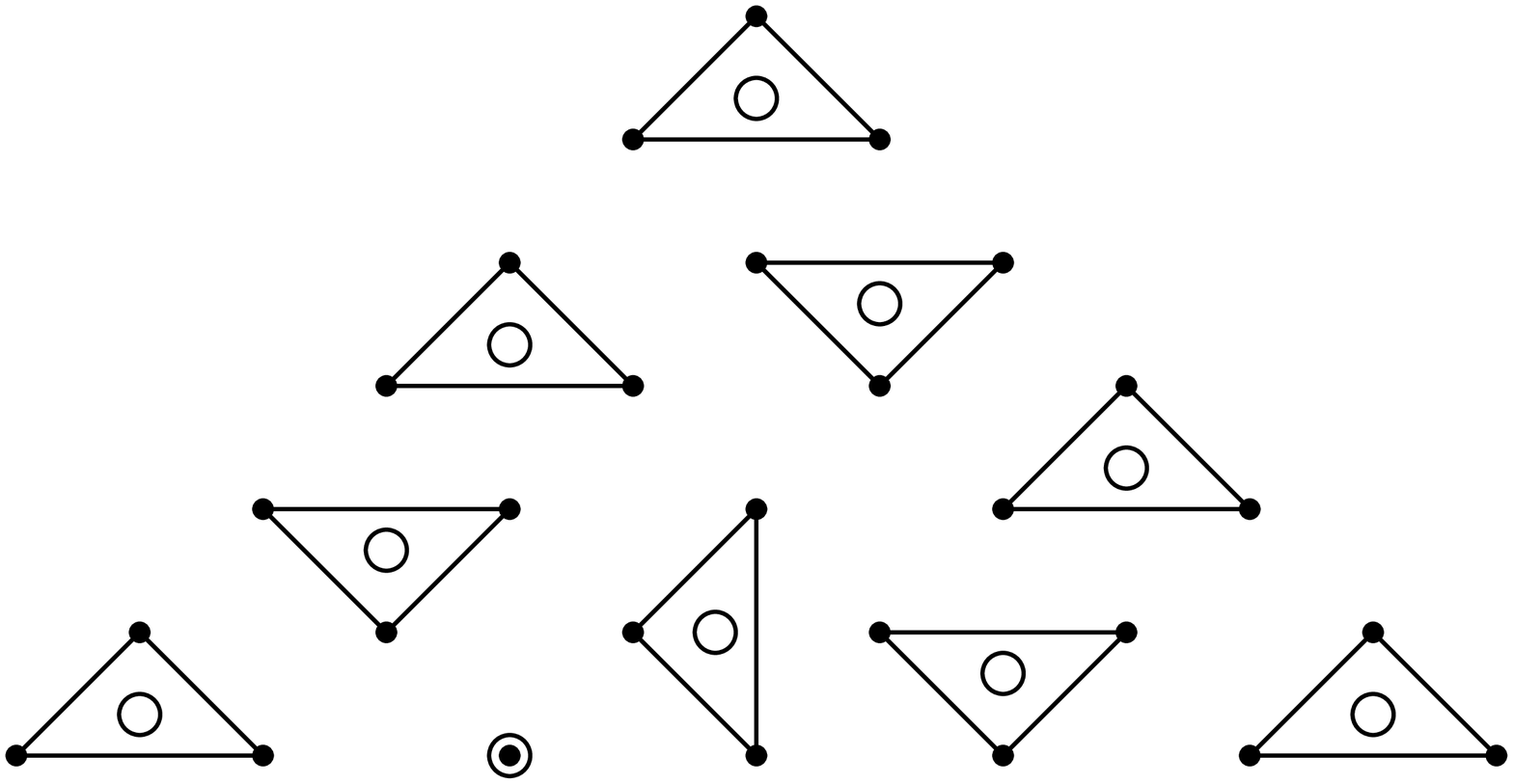}\label{fig:m3d6}}
\hfill
\subfigure[$d=7$]{\includegraphics[width=.4\textwidth]{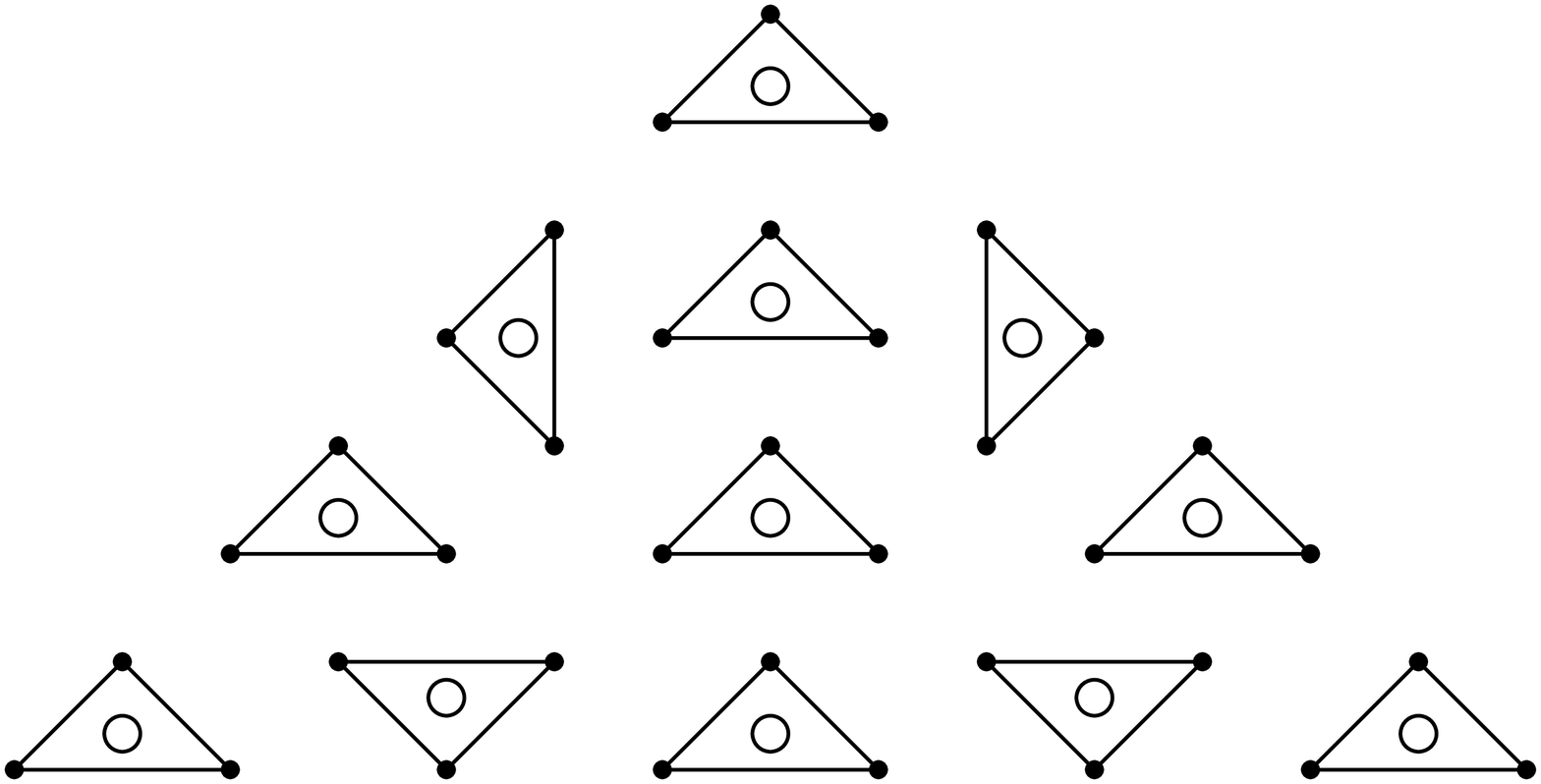}\label{fig:m3d7}}
\hfill
\subfigure[$d=8$]{\includegraphics[width=.4\textwidth]{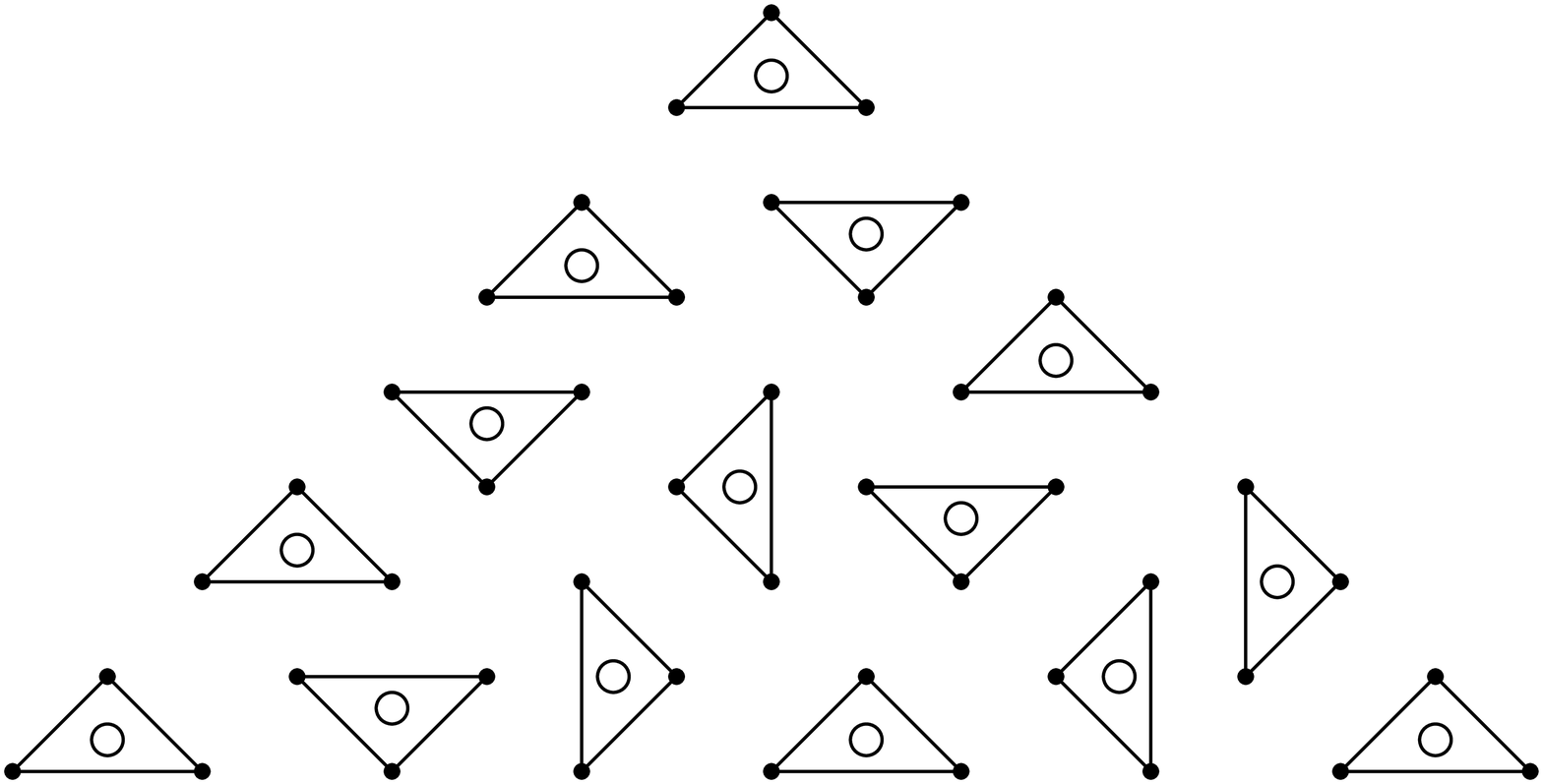}\label{fig:m3d8}}
\caption{Non-defective figures exist for $d \neq 2,4$.}
\label{fig:d1to8}
\end{figure} 

\begin{proof}
That $C$ is defective for $d=2,4$, is well known (for $d=4$ by the work of
Clebsch \cite{Clebsch1861}); in those two cases we will only show that the
defect is not more than $1$. We give a pictorial proof of the theorem:
the elements of $X$ form a triangle in the plane in $\RR^3$ where the
sum of the coordinates equals $d$. If we draw points $v_1,\ldots,v_k$
in that triangle, and if the points of $X$ lying in the Voronoi cell of
$v_i$ span an affine space of dimension $d_i \in \{0,1,2\}$, then $\sum_i
(d_i+1)$ is a lower bound on $\dim kC$ by Corollary \ref{cor:Secant}
and Lemma \ref{lm:Interrelations}. Moreover, we may and choose any
$2$-norm on the plane containing $X$. We choose to draw $X$ as in Figure
\ref{fig:d1to8}, and choose the $2$-norm for which circles really look
like circles in the plane. In this manner, Figure \ref{fig:m3d2} gives
lower bounds $3,5,6$ for the dimensions of $C,2C,3C$ when $d=2$: take
for the $v_i$ the midpoints of the open circles, and group the (black)
points of $X$ according to closest $v_i$; this results in a triangle,
an edge, and a single point.

\begin{figure}
\subfigure[From $d-2$ to $d \equiv 0 \mod 3$]{\includegraphics[width=.4\textwidth]{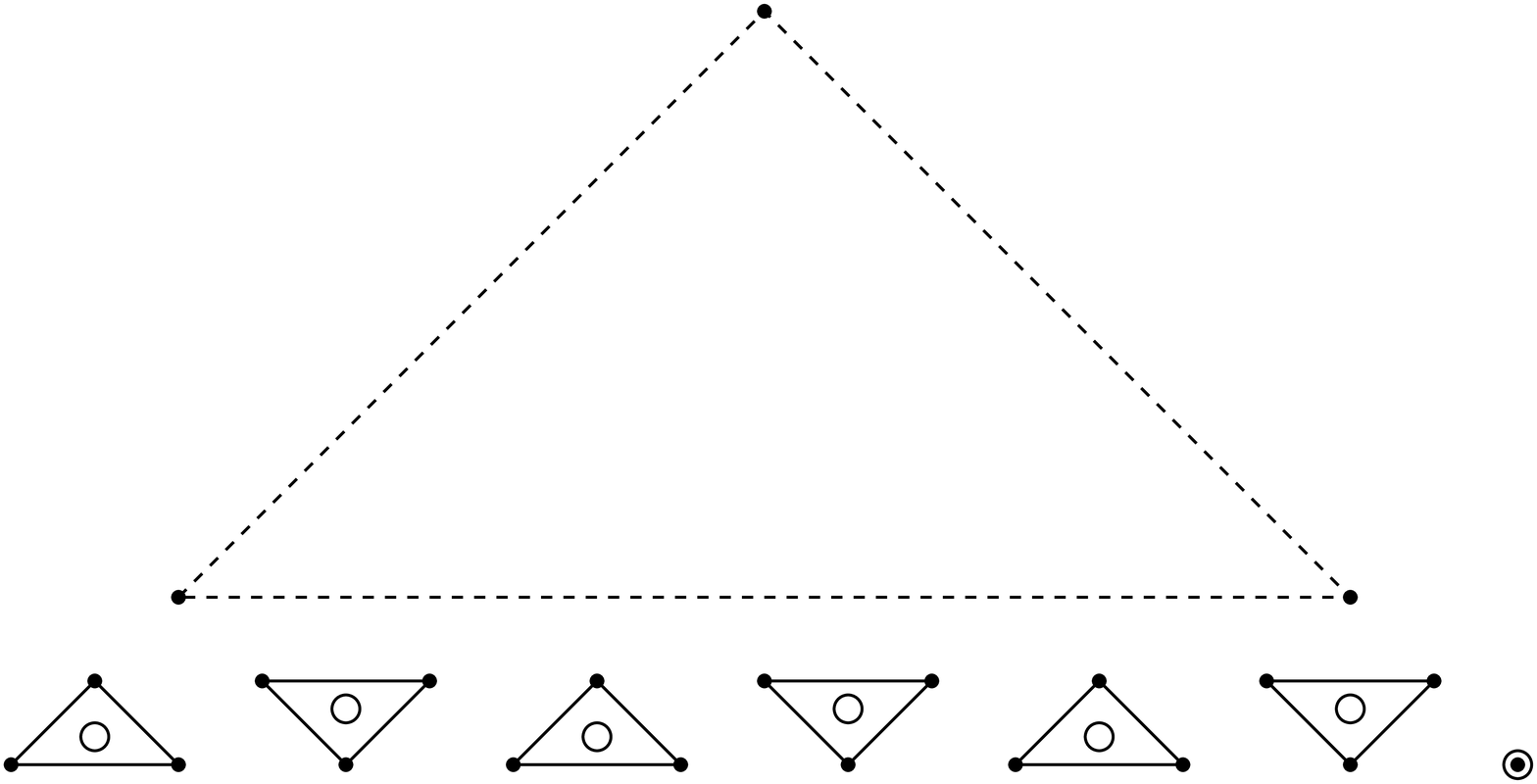}\label{fig:1naar0}}
\hfill
\subfigure[From $d-2$ to $d \equiv 1 \mod 3$]{\includegraphics[width=.4\textwidth]{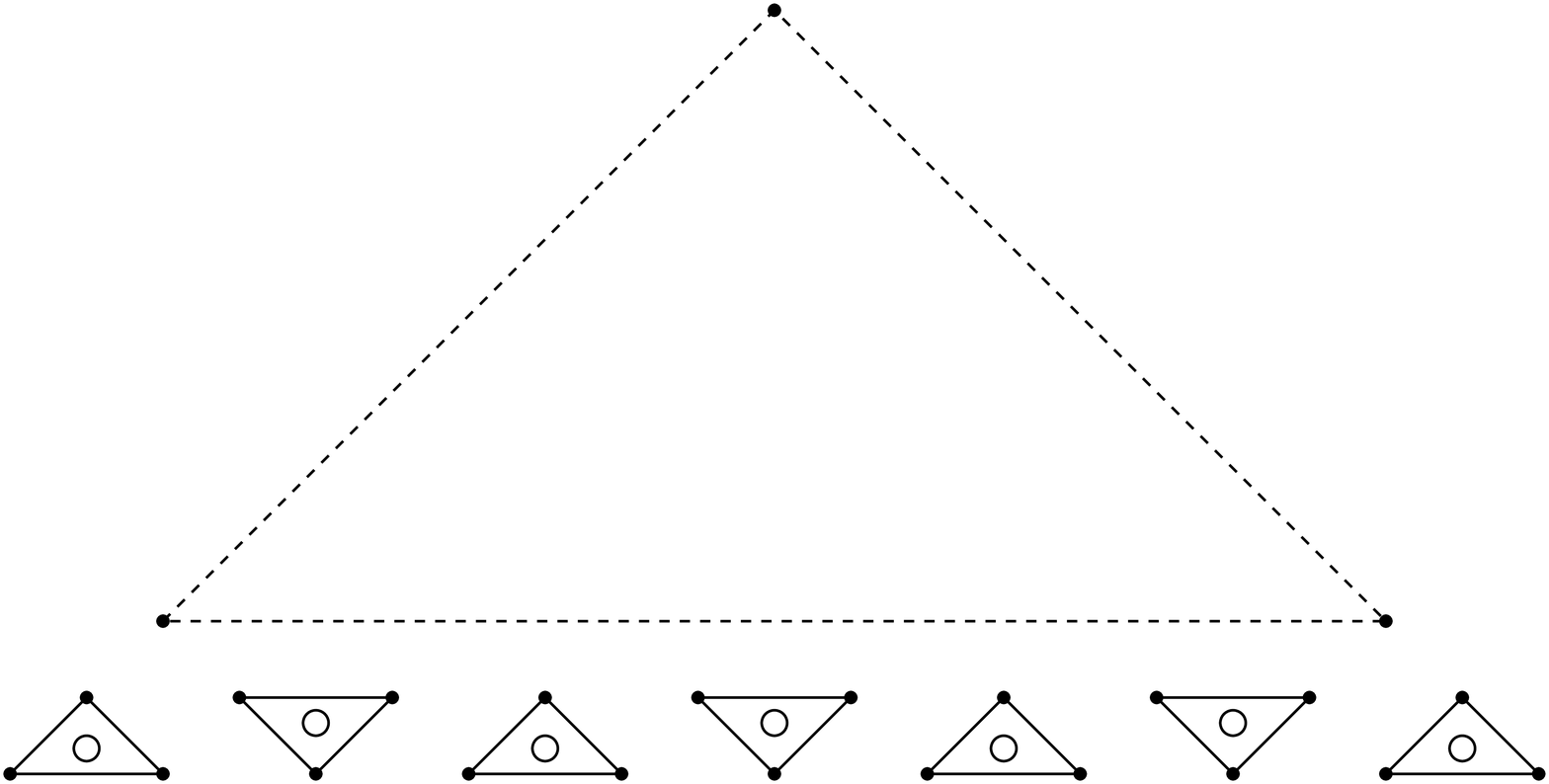}\label{fig:2naar1}}
\hfill
\subfigure[From $d-6$ to $d \equiv 2 \mod 3$]{\includegraphics[width=.4\textwidth]{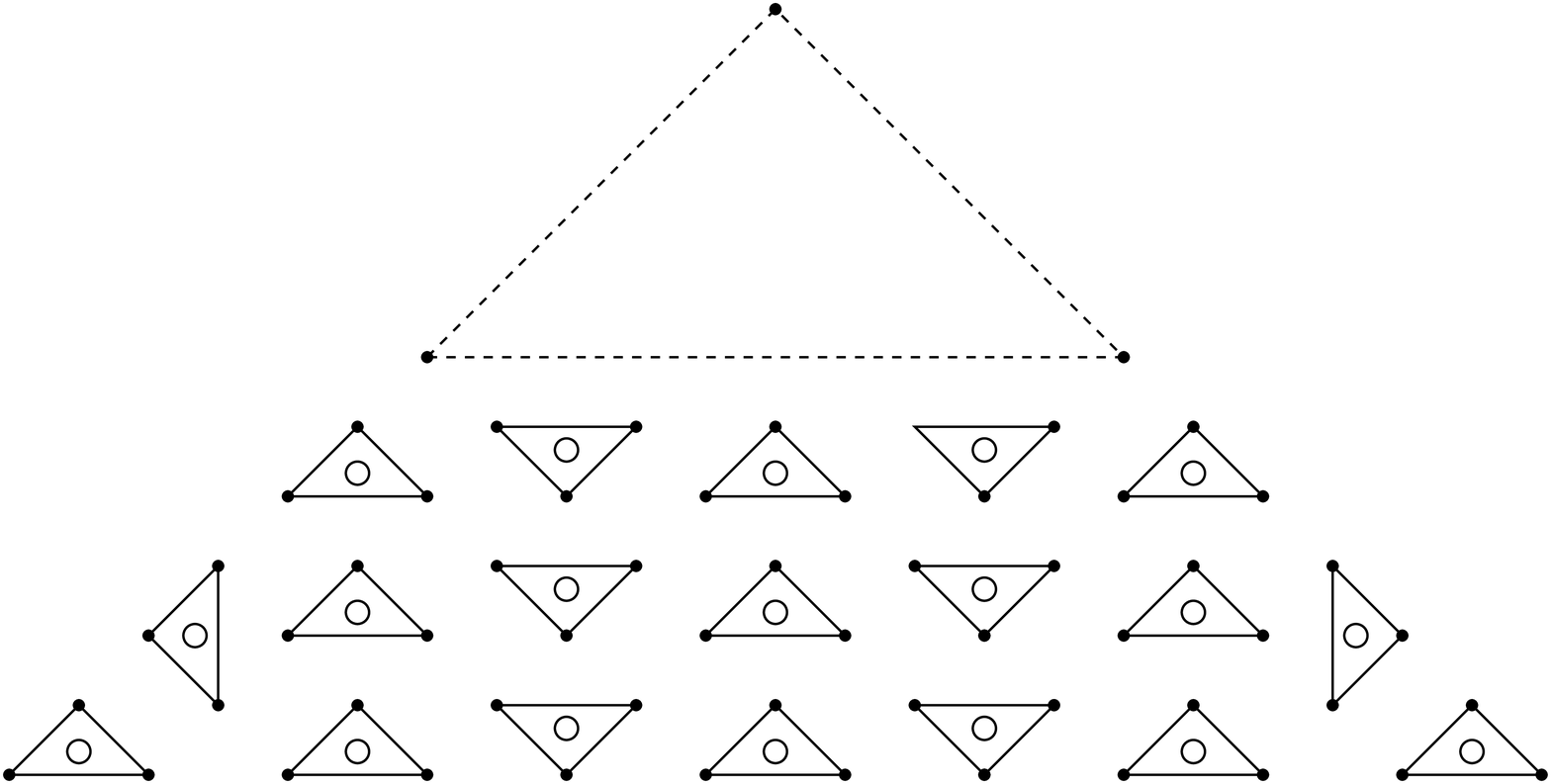}\label{fig:2naar2}}
\caption{Inductive construction of non-defective figures}
\end{figure}

We call such a picture non-defective if the groups are only
triangles, together with a single point if $d\equiv 0 \mod 3$; then
the picture shows that the corresponding $C$ is non-defective. Figures
\ref{fig:m3d3}--\ref{fig:m3d8} prove the theorem for $d=3,\ldots,8$;
all except Figure \ref{fig:m3d4} are non-defective. Now we proceed by
induction: we can produce a non-defective picture for $d \geq 9,\ d \equiv
0 \mod 3$ from a non-defective picture for $d-2$ as indicated in Figure
\ref{fig:1naar0}. Similarly, we construct a non-defective picture for
$d \geq 7,\ d \equiv 1 \mod 3$ from a non-defective picture for $d-2$
as in Figure \ref{fig:2naar1}. Finally, for $d \geq 11,\ d \equiv 2
\mod 3$, we construct a picture from a non-defective picture for $d-6$
as in Figure \ref{fig:2naar2}.  One readily verifies that this yields
non-defective figures for all $d \geq 8$ and hence proves the theorem.
\end{proof}

As promised, we compare the bound of Proposition \ref{prop:WeylOrbit}
to that of Corollary \ref{cor:Secant}.

\begin{prop} \label{prop:VeroneseWeyl}
The optimum of $\VP(X,k)$, relative
to the standard inner product, is at least that of Proposition
\ref{prop:WeylOrbit} applied to Veronese embeddings.
\end{prop}

\begin{proof}
Recall the bound of Proposition \ref{prop:WeylOrbit} for the cone of pure
powers: for indices $1\leq i_1<\ldots<i_k \leq n$ the lower bound on $\dim
kC$ in that proposition is the set of all words at $1$-distance at most
$2$ from the set $\{de_{i_1},\ldots,de_{i_k}\}=:B$. First suppose that $d
\geq 3$, so that the lower bound equals $km$. Choose $v_i:=de_{i_j}$,
set $v:=(v_1,\ldots,v_k) \in (\RR^m)^k$, and consider the Voronoi
diagram of the $v_i$. Then the $2$-distance between distinct $v_j$
is $d\sqrt{2}>2\sqrt{2}$, so that the Voronoi cell of $v_j$ contains
all $\alpha \in X$ at $2$-distance at most $\sqrt{2}$ from $v_j$. But
these are exactly the $m$ vectors in $X$ at $1$-distance at most $2$
from $de_{i_j}$, and their affine span $(m-1)$-dimensional; this shows
that $\sum_i (1+\dim \Aff_\RR D_i(v))$ is indeed at least the bound of
Proposition \ref{prop:WeylOrbit}.

Now suppose that $d=2$. First choose the $v_j$ as above.  Then $de_{i_j}$
still has $v_j$ as the unique closest point among all $v_l$, but $v_j$
and $v_l$ are equally close to $\alpha=e_{i_j}+e_{i_l}$. Now perturb all $v_j$s
slightly to $v_j':=v_j+\epsilon_j$ to resolve these ties, such that
relative to the Voronoi diagram of the $v_j'$ the vector $2e_{i_j}$
still lies in the cell of $v_j'$, while $e_{i_j}+e_{i_l}$ either lies
in the cell of $v_j'$ or in that of $v_l'$, for all $j \neq l$.  Then,
writing $v':=(v_j')_j$, each of the $\alpha$ at $2$-distance $\sqrt{2}$
from $B$ contributes exactly $1$ to the sum
\[ \sum_{i=1}^k (1+\dim \Aff_\RR D_i(v')), \] 
and therefore this sum is equal to $n+(n-1)+\ldots+(n-k+1)$, just like
the (exact) bound of Proposition \ref{prop:WeylOrbit}.
\end{proof}

The proof of Theorem \ref{thm:Veronesem3} shows that the optimum of $\VP$
is in fact much stronger than that of Proposition
\ref{prop:WeylOrbit}: while the latter bound only takes corners
of $X$ into account, the former bound also exploits the interior
of $X$. We will encounter a similar phenomenon with Segre varieties.

\subsection{Segre varieties}
Let $m,d$ be natural numbers and set $V:=(K^m)^{\otimes d}$. Let $X$ be
the set of $d \times m$-matrices of natural numbers whose row sums are
all $1$. For $\alpha \in X$ we set $e_\alpha:=e_{j_1}\otimes\ldots\otimes
e_{j_d}$, where $j_i$ is the column containing the $1$ in the $i$-th
row; the $e_\alpha$ form our basis of $V$. Let $C$ be the closed cone
of of pure tensors in $V$. Writing $M_{d,m}(K)$ for the space
of $d\times m$-matrices with entries in $K$, we parametrise $C$ by the
polynomial map $M_{d,m}(K) \rightarrow V$ sending $x$ to
\[ (x_{11}e_1+\ldots+x_{1m}e_m)\otimes \ldots \otimes
	(x_{d1}e_1+\ldots+x_{dm}e_m) = \sum_{\alpha \in X} x^\alpha
	e_\alpha. \]
By Corollary \ref{cor:Secant}, the optimum of $\LP(\{\{\alpha\}\}_{\alpha
\in X}, k)$ is a lower bound on $\dim kC$ for all $k$.  Moreover, $X$
is contained in the affine hyperspace of $M_{d,m}(\RR)$ where all row
sums are $1$, so that also the optima of $\AP(\{\{\alpha\}\}_{\alpha
\in X}, k)$ and $\VP(X, k)$ are lower bounds on $\dim kC$ by Lemma
\ref{lm:Interrelations}. This last statement is the content of Proposition
\ref{prop:MultOne} for the minimal orbit of the representation of
$\lieg{GL}_m^d$ on $V$.  First we show that the last lower bound is at
least as good as that of Proposition \ref{prop:WeylOrbit}.

\begin{prop} \label{prop:SegreWeyl}
The optimum of $\VP(X,k)$, relative to the
standard inner product on $M_{d,r}(\RR)$, is at least the lower bound
of Proposition \ref{prop:WeylOrbit} applied to Segre powers.
\end{prop}

\begin{proof}
The proof is very similar to that of Proposition
\ref{prop:VeroneseWeyl}: the bijection $\phi$ from $\{1,\ldots,m\}^d$
to $X$ sending $(j_1,\ldots,j_d)$ to the matrix having $1$s on
the positions $(i,j_i)$ and $0$s elsewhere has the property that
$||\phi(w_1)-\phi(w_2)||_2=\sqrt{2b}$ if $b$ is the Hamming distance
between $w_1$ and $w_2$. Let $w_1,\ldots,w_k$ be distinct $m$-ary words
of length $d$, set $v_i:=\phi(w_i)$ and $v:=(v_1,\ldots,v_k)$. The lower
bound on $\dim kC$ of Proposition \ref{prop:WeylOrbit} is then size of
the set
\[ S:=\{\alpha \in X \mid ||\alpha-v_i|| \leq \sqrt{2}
	\text{ for some } i\}. \]
Note that for every $i$ the elements of $S$ at $2$-distance $\leq
\sqrt{2}$ from $v_i$ form an affinely independent set: they are either
equal to $v_i$ or obtained from $v_i$ by moving a $1$ within its row.
Some elements of $S$ may have distance $\sqrt{2}$ to two distinct $v_i$.
Perturbing $v$ slightly to resolve these ties yields a $v'$ such that
$S \subseteq \bigcup_i D_i(v')$ and every $S \cap D_i(v')$ is contained
in the closed $2$-ball of radius $\sqrt{2}$ from $v_i$---hence affinely
independent. Hence
$\sum_{i=1}^k (1+\dim \Aff_\RR D_i(v'))$
is at least $|S|$, which is the lower bound from Proposition
\ref{prop:WeylOrbit}.
\end{proof}

The following proposition shows that for Segre powers, too, the optimum
of $\VP$ can be strictly larger than the lower bound of Proposition
\ref{prop:WeylOrbit}.

\begin{prop}
The $6$-th Segre power of the projective line is non-defective.
\end{prop}

In \cite{Catalisano05b} it is proved that all higher secant varieties of
Segre powers of the projective line are non-defective, except possibly
for one higher secant variety of each Segre power. The statement
of the proposition is the first case not covered by the theorem in
\cite{Catalisano05b}, and what follows is the first computer-free proof
of that statement.

It will be convenient to work with the set $Y \subseteq M_{d,m-1}(\RR)$
obtained from $X$ by deleting the last column of every element of $X$;
note that this operation defines an affine equivalence from $X$ to $Y$,
so that it does not affect the optimum of $\AP$. We will work with $\VP$
relative to the standard inner product on $M_{d,m-1}(\RR)$; note that
its restriction to $Y$ is {\em not} equal to that on $M_{d,m}(\RR)$
restricted to $X$ and transferred to $Y$ by the affine equivalence---but
in $\VP$ we are free to choose our inner product!

\begin{proof}
We are in the situation where $m=2$ and $d=6$. We have $\dim C=7$
and $\dim V=64$, so we have to show that $9C$ has the expected
dimension $63$. This is impossible using only the rook covering
bound of Proposition \ref{prop:WeylOrbit}, because the maximal size of
a binary code of length $6$ and Hamming distance $3$ is $8$ (see, e.g.,
\cite{Best78}). However, we will use such a code, and then complement it
with a further point to take care of the points outside the Hamming
balls of radius $1$ around its codewords.

More specifically, note that $Y=\{0,1\}^6$. Let $B$ be the set of all
vectors $v$ in $Y$ with
\[ Hv=0 \mod 2,\quad \text{ where } 
	H=\begin{bmatrix}
		0&0&0&1&1&1\\
		0&1&1&0&0&1\\
		1&0&1&0&1&0\\
	\end{bmatrix}
\]
is the {\em parity check matrix} of $B$. From the fact that the columns
of $H$ are distinct modulo $2$ one readily concludes that the minimal
Hamming distance between elements of $B$ is $3$. Explicitly, $B$ consists
of the rows $v_1,\ldots,v_8$ of the following matrix:
\begin{align*} 
	\begin{bmatrix}
		0&0&0&0&0&0\\
		1&0&0&1&1&0\\
		0&1&0&1&0&1\\
		0&0&1&0&1&1\\
		1&1&0&0&1&1\\
		1&0&1&1&0&1\\
		0&1&1&1&1&0\\
		1&1&1&0&0&0
	\end{bmatrix}
\end{align*}
Now set $v_9:=(\frac12,\ldots,\frac12)$ and consider the Voronoi
diagram of $v=(v_1,\ldots,v_9)$. Any element $w$ of $Y$ has $2$-distance
$\sqrt{3/2}$ to $v_9$. If $w$ has Hamming distance $1$ to some $v_i$
with $i<9$, then its $2$-distance to $v_i$ is also $1$, and $w$ lies
in the Voronoi cell of $v_i$. Otherwise, $w$ has $2$-distance at least
$\sqrt{2}$ to every $v_i$ with $i<9$, hence $w$ lies in the Voronoi
cell of $v_9$. We have to check that each $D_i(v)$ spans an affine
space of dimension $6$. For $i<9$ this is clear, because $D_i(v)$
is affinely equivalent to the Hamming ball of radius $1$ around $0$,
which apart from $0$ contains all standard basis vectors. Finally,
$D_9(v)$ contains $64-8*7=8$ words, and it is easy to see that these are
precisely the words of the form $(1,1,1,1,1,1)-v_i$ with $i<9$: indeed,
these $8$ words give syndrome $(1,1,1)^t$ when multiplied from the right
by $H$ and taken modulo $2$, so they are not at Hamming distance $1$
from any element of $B$. Hence $D_i(v)$ is affinely equivalent to $B$,
and a direct computation shows that the affine span of $B$ is the full
space $\RR^6$. We conclude that $9C$ is, indeed, non-defective.
\end{proof}

\subsection{Grassmannians} \label{ssec:Grassmannians}

Let $m,d$ be natural numbers with $d \leq m/2$ and set $V:=\Wedge^d
(K^m)$. Let $X$ be the set of all $d$-subsets of $\{1,\ldots,m\}$ and
for $J \in X$ with elements $i_1<\ldots<i_d$ set $e_J:= e_{i_1} \wedge
\ldots \wedge e_{i_d}$; the $e_J$ form our basis of $V$. Let $C$ be the
closed cone which is the image of $M_{d,m}(K)$ under the map sending $x$ to
\[ (x_{11} e_1 + \ldots + x_{1m} e_m) \wedge \ldots \wedge
(x_{d1} e_1 + \ldots + x_{dm} e_m) =\sum_J \det(x_J) e_J; \]
here $x_J$ is the $d \times d$-matrix obtained from $x$ by taking only
the columns corresponding to elements of $J$. The exponent vectors of the
monomials appearing in the coefficient $\det(x_J)$ of $e_J$ are precisely
those matrices in $M_{d,m}(\NN)$ that are $0$ outside the columns in $J$
and whose restriction to the columns in $J$ are permutation matrices;
call this collection $A_J$. Then Corollary \ref{cor:Secant} says that
the optimum of $\LP((A_J)_J,k)$ is a lower bound on $\dim kC$. Again,
all $A_J$ lie in a common affine space not through zero, so that
also $\AP$ yields a lower bound. Since the $A_J$ are not singletons,
we cannot immediately use $\VP$. In fact, one can still prove that
the optimum of $\AP((A_J)_J,k)$ is greater than the lower bound of
Proposition \ref{prop:WeylOrbit}, but it is slightly more involved. As
it fits better with forthcoming work dealing with more general secant
dimensions of minimal orbits, we omit it here.

\section{Conclusion} \label{sec:Conclusion}
The tropical approach to secant dimensions shows very promising
results when tested on concrete minimal orbits, especially those in
representations where all weight spaces are one-dimensional---like the
Veronese, Segre, or Pl\"ucker embeddings. The approach leads to exciting
combinatorial-polyhedral questions. In particular, the approach yields a
nice pictorial proof of the non-defectiveness of most Veronese surfaces,
one of the main results of \cite{Hirschowitz85}; a similar proof for
general Veronese embeddings would be an attractive alternative to
\cite{Alexander95}.

Still, these polyhedral-combinatorial questions are mostly open, and
there is a lot of space for further research. In particular:

\begin{enumerate}

\item Using Terracini's lemma, one can compute $\dim kC$ by computing
the rank of the addition map $C^k \rightarrow kC$ at a {\em generic}
point. However, for this one has to compute the rank of a large matrix
(of size the dimension of the representation); see \cite{Baur06}. On
the other hand, the tropical approach only needs the ranks of several
smaller matrices (of size the dimension of $C$), but it only works at a
carefully selected point where the rank of the differential of a tropical
polynomial map is maximal. Is there a method in between, which does work
at random points but only requires ranks of small matrices?

\item The tropical approach depends on the chosen bases: to prove anything
substantial, it seems wise to choose nice bases of both the representation
and the parametrising space. For minimal orbits in representations where
not all weight spaces are one-dimensional, it is unclear which bases of
the representation one should use.  The question of whether there exist
bases for which the tropical method works well is a very exciting one!

\item Proposition \ref{prop:MultOne} gives very interesting lower bounds
for secant dimensions of minimal orbits where the root system of the
negative unipotent radical is linearly independent. In particular,
these bounds apply to Segre products of Veronese embeddings. I do not
know of an example where the lower bound is not sharp.

\end{enumerate}

%\bibliographystyle{plain}
%\bibliography{diffeq,jan}

\begin{thebibliography}{10}

\bibitem{Alexander88}
J.~Alexander.
\newblock Singularit\'{e}s imposables en position g\'en\'erale \`a une
  hypersurface projective.
\newblock {\em Compos. Math.}, 68(3):305--354, 1988.

\bibitem{Alexander92}
J.~Alexander and A.~Hirschowitz.
\newblock La m\'ethode d'{H}orace \'eclat\'ee: Application \`a l'interpolation
  en degr\'e quatre.
\newblock {\em Invent. Math.}, 107(3):585--602, 1992.

\bibitem{Alexander95}
J.~Alexander and A.~Hirschowitz.
\newblock Polynomial interpolation in several variables.
\newblock {\em J. Algebr. Geom.}, 4(2):201--222, 1995.

\bibitem{Baur06}
Karin Baur, Willem de~Graaf, and Jan Draisma.
\newblock Secant dimensions: computations and conjectures.
\newblock 2006.
\newblock Preprint.

\bibitem{Baur04}
Karin Baur and Jan Draisma.
\newblock Higher secant varieties of the minimal adjoint orbit.
\newblock {\em J. Algebra}, 280:743--761, 2004.

\bibitem{Best78}
M.R. Best, A.E. Brouwer, F.Jessie MacWilliams, Andrew~M. Odlyzko, and Neil~J.A.
  Sloane.
\newblock Bounds for binary codes of length less than 25.
\newblock {\em IEEE Trans. Inf. Theory}, 24:81--93, 1978.

\bibitem{Bieri84}
Robert Bieri and J.R.J. Groves.
\newblock The geometry of the set of characters induced by valuations.
\newblock {\em J. Reine Angew. Math.}, 347:168--195, 1984.

\bibitem{Borel91}
Armand Borel.
\newblock {\em Linear Algebraic Groups}.
\newblock Springer-Verlag, New York, 1991.

\bibitem{Bosch84}
S.~Bosch, U.~G\"{u}ntzer, and R.~Remmert.
\newblock {\em Non-Archimedean Analysis}, volume 261 of {\em Grundlehren der
  mathematischen Wissenschaften}.
\newblock Springer-Verlage, Berlin, Heidelberg, 2004.

\bibitem{Catalisano2002}
M.V. Catalisano, A.V. Geramita, and A.~Gimigliano.
\newblock {On the rank of tensors, via secant varieties and fat points.}
\newblock In {\em {Geramita, A. V. (ed.), Zero-dimensional schemes and
  applications. Proceedings of the workshop, Naples, Italy, February 9--12,
  2000. Kingston: Queen's University.}}, volume {123} of {\em {Queen's Pap.
  Pure Appl. Math.}}, pages {135--147}, 2002.

\bibitem{Catalisano05b}
M.V. Catalisano, A.V. Geramita, and A.~Gimigliano.
\newblock Higher secant varieties of the {S}egre varieties
  $\mathbb{P}^1\times\dots\times\mathbb{P}^1$.
\newblock {\em J. Pure Appl. Algebra}, 201(1--3):367--380, 2005.

\bibitem{Catalisano05}
M.V. Catalisano, A.V. Geramita, and A.~Gimigliano.
\newblock Secant varieties of grassmann varieties.
\newblock {\em Proc. Am. Math. Soc.}, 133(3):633--642, 2005.

\bibitem{Clebsch1861}
A.~Clebsch.
\newblock {\"U}ber {K}urven vierter {O}rdnung.
\newblock {\em J. Reine Angew. Math.}, 59:125--145, 1861.

\bibitem{Develin06}
Mike Develin.
\newblock Tropical secant varieties of linear spaces.
\newblock {\em Discrete Comput. Geom.}, 35(1):117--129, 2006.

\bibitem{Ehrenborg00}
R.~Ehrenborg.
\newblock On apolarity and generic canonical forms.
\newblock {\em J. of Algebra}, 213:167--194, 1999.

\bibitem{Einsiedler04}
Manfred Einsiedler, Mikhail Kapranov, and Douglas Lind.
\newblock Non-archimedean amoebas and tropical varieties.
\newblock 2004.
\newblock Preprint, available from \verb+http://arxiv.org/abs/math.AG/0408311+.

\bibitem{Fresnel04}
Jean Fresnel and Marius van~der Put.
\newblock {\em Rigid analytic geometry and its applications}, volume 218 of
  {\em Progress in Mathematics}.
\newblock Birkh\"auser, Boston, 2004.

\bibitem{Hirschowitz85}
Andr\'e Hirschowitz.
\newblock La methode d'{H}orace pour l'interpolation \`a plusieurs variables.
\newblock {\em Manuscr. Math.}, 50:337--388, 1985.

\bibitem{Humphreys72}
James~E. Humphreys.
\newblock {\em Introduction to {L}ie algebras and representation theory}.
\newblock Springer, Berlin, 1972.

\bibitem{Humphreys81}
James~E. Humphreys.
\newblock {\em Linear algebraic groups}, volume~21 of {\em Graduate Texts in
  Mathematics}.
\newblock Springer-Verlag, New York - Heidelberg - Berlin, 1981.

\bibitem{Lint99}
J.H.~van Lint.
\newblock {\em Introduction to coding theory}, volume~86 of {\em Graduate Texts
  in Mathematics}.
\newblock Springer, Berlin, 1991.

\bibitem{Mikhalkin04}
Grigory Mikhalkin.
\newblock Amoebas of algebraic varieties and tropical geometry.
\newblock In Simon~Donaldson {\em et al.}, editor, {\em Different faces of
  geometry}, volume~3 of {\em International Mathematical Series}, pages
  257--300, New York, 2004. Kluwer Academic/Plenum Publishers.

\bibitem{Mikhalkin06}
Grigory Mikhalkin.
\newblock Tropical geometry and its applications.
\newblock 2006.
\newblock Preprint, available from
  \verb+http://www.arxiv.org/abs/math.AG/0601041+.

\bibitem{Pachter05}
Lior Pachter and Bernd Sturmfels, editors.
\newblock {\em Algebraic statistics for computational biology}, Cambridge,
  2005. Cambridge University Press.

\bibitem{RichterGebert05}
J\"urgen Richter-Gebert, Bernd Sturmfels, and Thorsten Theobald.
\newblock First steps in tropical geometry.
\newblock In G.~L.~{\em et al.} Litvinov, editor, {\em Idempotent mathematics
  and mathematical physics. Proceedings of the international workshop, Vienna,
  Austria, February 3-10, 2003}, volume 377 of {\em Contemporary Mathematics},
  pages 289--317, Providence, RI, 2005. AMS.

\bibitem{Speyer04}
David Speyer and Bernd Sturmfels.
\newblock The tropical grassmannian.
\newblock {\em Adv. Geom.}, 4(3):389--411, 2004.

\bibitem{Sturmfels02}
Bernd Sturmfels.
\newblock {\em Solving systems of polynomial equations}, volume~97 of {\em CBMS
  Regional Conference Series in Mathematics}.
\newblock AMS, Providence, RI, 2002.

\bibitem{Sturmfels05}
Bernd Sturmfels and Seth Sullivant.
\newblock Combinatorial secant varieties.
\newblock {\em Quart. J. Pure Applied Math.}, 2005.
\newblock To appear; preprint available from
  \verb+http://arxiv.org/abs/math.AC/0506223+.

\bibitem{Terracini11}
A.~Terracini.
\newblock {Sulle $V_k$ per cui la varietà degli $S_h(h+1)$-seganti ha
  dimensione minore dell' ordinario.}
\newblock {\em Palermo Rend.}, 31:392--396, 1911.

\bibitem{Zak93}
F.L. Zak.
\newblock {\em {Tangents and secants of algebraic varieties.}}, volume 127 of
  {\em Translations of Mathematical Monographs}.
\newblock American Mathematical Society (AMS), Providence, RI, 1993.

\end{thebibliography}

\end{document}